\documentclass[12pt]{amsart}
\usepackage{amsmath,amscd,amssymb,amsfonts}
\newcommand{\ver}{June 12, 2007, v.4}
\newcommand{\scirc}{\raise.2ex\hbox{${\scriptstyle\circ}$}}
\newcommand{\ssbull}{\raise.2ex\hbox{${\scriptscriptstyle\bullet}$}}
\newcommand{\moplus}{\hbox{$\bigoplus$}}
\newcommand{\mprod}{\hbox{$\prod$}}
\newcommand{\bD}{{\mathbb D}}
\newcommand{\bN}{{\mathbb N}}
\newcommand{\bP}{{\mathbb P}}
\newcommand{\bQ}{{\mathbb Q}}
\newcommand{\bZ}{{\mathbb Z}}
\newcommand{\boL}{{\mathbf L}}
\newcommand{\boR}{{\mathbf R}}
\newcommand{\cC}{{\mathcal C}}
\newcommand{\cE}{{\mathcal E}}
\newcommand{\cF}{{\mathcal F}}
\newcommand{\cH}{{\mathcal H}}
\newcommand{\cK}{{\mathcal K}}
\newcommand{\cL}{{\mathcal L}}
\newcommand{\cS}{{\mathcal S}}
\newcommand{\of}{\overline{f}}
\newcommand{\ok}{\overline{k}}
\newcommand{\oK}{\overline{K}}
\newcommand{\te}{\widetilde{e}}
\newcommand{\ti}{\widetilde{i}}
\newcommand{\tK}{\widetilde{K}}
\newcommand{\tX}{\widetilde{X}}
\newcommand{\rank}{\hbox{{\rm rank}}}
\newcommand{\can}{\hbox{{\rm can}}}
\newcommand{\Var}{\hbox{{\rm Var}}}
\newcommand{\IC}{\hbox{{\rm IC}}}
\newcommand{\CH}{\hbox{{\rm CH}}}
\newcommand{\Spec}{\hbox{{\rm Spec}}\,}
\newcommand{\Perv}{\hbox{{\rm Perv}}}
\newcommand{\Ext}{\hbox{{\rm Ext}}}
\newcommand{\Gr}{\text{{\rm Gr}}}
\newcommand{\Imm}{\hbox{{\rm Im}}}
\newcommand{\Ker}{\hbox{{\rm Ker}}}
\newcommand{\Coker}{\hbox{{\rm Coker}}}
\newcommand{\Hom}{\hbox{{\rm Hom}}}
\newcommand{\cHom}{{\mathcal H}om}
\newcommand{\Aut}{\hbox{{\rm Aut}}}

\newcommand{\Gal}{\hbox{{\rm Gal}}}
\newcommand{\Sh}{\hbox{{\rm Sh}}}
\newcommand{\TC}{\hbox{{\rm TC}}}
\newcommand{\supp}{\hbox{{\rm supp}}\,}

\newcommand{\red}{\text{\rm red}}
\newcommand{\sm}{\text{\rm sm}}
\newcommand{\gnr}{\text{\rm gnr}}
\newcommand{\gen}{\text{\rm gen}}
\newcommand{\inv}{\text{\rm inv}}
\newcommand{\van}{\text{\rm van}}
\newcommand{\homm}{\text{\rm hom}}
\newcommand{\simto}{\buildrel\sim\over\longrightarrow}

\begin{document}
\title{Tate Conjecture and Mixed Perverse Sheaves}
\author{Morihiko Saito}
\address{RIMS Kyoto University, Kyoto 606-8502 Japan}
\email{msaito@kurims.kyoto-u.ac.jp}
\date{\ver}
\begin{abstract}
Using the theory of mixed perverse sheaves, we extend arguments on
the Hodge conjecture initiated by Lefschetz and Griffiths to the
case of the Tate conjecture,
and show that the Tate conjecture for divisors is closely related
to the de Rham conjecture for nonproper varieties, finiteness of
the Tate-Shafarevich groups, and also to some conjectures in the
analytic number theory.
\end{abstract}
\maketitle

\bigskip
\centerline{Dedicated to John Tate}

\bigskip\bigskip
\centerline{\bf Introduction}

\bigskip\noindent
Let
$ k $ be a finitely generated field over
$ \bQ $,
$ \ok $ an algebraic closure of
$ k $,
and
$ G_{k} = \Gal(\ok/k) $.
Let
$ X $ be a smooth projective variety of pure dimension
$ n $ over
$ k $, and
$ X_{\ok} = X\otimes_{k}\ok $.
We will denote by
$ \TC(X/k,p) $ the Tate conjecture ([35], [37]) which states the
surjectivity of the cycle map
$$
cl : \CH^{p}(X)\otimes \bQ_{l}\to
H^{2p}(X_{\ok}, \bQ_{l}(p))^{G_{k}}.
\leqno(0.1)
$$
Here
$ \CH^{p}(X) $ is the Chow group of algebraic cycles of codimension
$ p $ on
$ X $,
and the right-hand side is the invariant part of the (Tate twisted)
\'etale cohomology group by the action of
$ G_{k} $.
Let
$ D $ be a smooth divisor on
$ X $,
and set
$$
\aligned
H^{j}(D_{\ok}, \bQ_{l})_{X}
&= \Coker(H^{j}(X_{\ok},
\bQ_{l})\to H^{j}(D_{\ok}, \bQ_{l})),
\\
H^{j}(X_{\ok}, \bQ_{l})^{D}
&= \Ker(H^{j}(X_{\ok}, \bQ_{l})\to H^{j}(D_{\ok},
\bQ_{l})),
\endaligned
$$
so that we have an exact sequence compatible with the Galois action
$$
0 \to H^{2p-1}(D_{\ok}, \bQ_{l})_{X} \to
{H}_{c}^{2p}(X_{\ok} \setminus D_{\ok},
\bQ_{l})\to H^{2p}(X_{\ok}, \bQ_{l})^{D} \to
0.
\leqno(0.2)
$$

For
$ c \in (H^{2p}(X_{\ok}, \bQ_{l}(p))^{D})^{G_{k}} \,(=
\Hom_{G_{k}}(\bQ_{l}, H^{2p}(X_{\ok}, \bQ_{l}(p))^{D})) $,
we denote by
$$
e(c) \in \Ext^{1}(\bQ_{l}, H^{2p-1}(D_{\ok},
\bQ_{l}(p))_{X})
$$
the extension class in the category of
$ \bQ_{l} $-modules with action of
$ G_{k} $,
which is obtained by taking the pull-back of (0.2) by
$ c $
(see [26], [28] for the Hodge case).
Let
$$
\CH_{\homm}^{p}(D) = \Ker(cl : \CH^{p}(D) \to H^{2p}(D_{\ok},
\bQ_{l}(p))).
$$
By a similar argument, we get the Abel-Jacobi map
$$
\CH_{\homm}^{p}(D) \to
\Ext^{1}(\bQ_{l}, H^{2p-1}(D_{\ok}, \bQ_{l}(p))),
$$
which has been obtained in [19], see also (2.1) below.
In the Hodge setting, this construction is essentially due to Deligne,
see [14] (and also [26], [28]).
Taking the composition with the canonical projection we get
$$
\CH_{\homm}^{p}(D)\otimes \bQ_{l}\to
\Ext^{1}(\bQ_{l}, H^{2p-1}(D_{\ok}, \bQ_{l}(p))_{X}),
\leqno(0.3)
$$
where the extension group may be replaced with Galois cohomology.
We have the following

\medskip\noindent
{\bf 0.4.~Conjecture.} For any
$ c \in (H^{2p}(X_{\ok}, \bQ_{l}(p))^{D})^{G_{k}} $,
the above extension class
$ e(c) $ belongs to the image of (0.3).

\medskip

Note that (0.4) follows from
$ \TC(X/k,p) $,
because
$ e(c) $ coincides with the image of the restriction of
$ \zeta $ to
$ D $ by (0.3) if
$ c $ is the cycle class of a cycle
$ \zeta $ on
$ X $, see [28], 1.8.
Conversely, we can reduce the Tate conjecture to (0.4) using a Lefschetz
pencil
(by induction on
$ \dim X) $,
where
$ D $ is the generic fiber of the Lefschetz pencil and the base field
$ k $ is replaced by the rational function field
$ k(t) $ of one variable.
Indeed,
$ X $ can be embedded in a projective space, and we have a Lefschetz
pencil
$ f : \tX \to S = \bP^{1} $, where
$ \pi : \tX \to X $ is the blow up along the intersection
$ A $ of two general hyperplane sections defined over
$ k $.
Then
$ \pi \times f : \tX \to X\times_{k}S $ is a closed
embedding
so that the closed fibers of
$ \of = f\otimes_{k}\ok : \tX_{\ok} \to
S_{\ok} $ are identified with the hyperplane sections of
$ X_{\ok} $ containing
$ A_{\ok} $.
For a closed point
$ s $ of
$ S $,
we will denote by
$ X_{s} $ the fiber of
$ f $ over
$ s $.
The generic fiber of
$ f $ will be denoted by
$ Y $.
It is smooth projective over
$ \Spec K $ with
$ K := k(S) = k(t) $.
Let
$ X_{K} = X\otimes_{k}K $ so that
$ Y $ is a closed subvariety of
$ X_{K} $.
Let
$ U $ be a nonempty open subvariety of
$ S $ on which
$ f $ is smooth, and
$ |U| $ the set of closed points of
$ U $.
Then we can prove (see (2.5)):

\medskip\noindent
{\bf 0.5.~Lemma.} {\it
{\rm (i)} In case
$ n < 2p $,
$ \TC(X/k,p) $ is true if
$ \TC(X_{s}/k(s), p-1) $ is true for some
$ s \in |U| $.

\medskip\noindent
{\rm (ii)} In case
$ n > 2p $,
$ \TC(X/k,p) $ is true if
$ \TC(Y/K,p) $ and
$ \TC(X_{s}/k(s), p-1) $ are true for some
$ s \in |U| $.
}

\medskip
So the Tate conjecture is reduced to the case
$ n = 2p $ by induction.
Assume the projective embedding of
$ X $ is sufficiently ample so that we have by [20], 6.3, 6.4
$$
H^{n-1}(Y_{\oK}, \bQ_{l})_{X} \ne 0,
\leqno(0.6)
$$
which is equivalent to the constantness of
$ R^{j}\of_{*}\bQ_{l} $ on
$ S_{\ok} $ for
$ j \ne n - 1 $, see (2.4) below.
We have the following (see (2.6) below):

\medskip\noindent
{\bf 0.7.~Theorem.} {\it
For
$ n = 2p $,
$ \TC(X/k,p) $ is true, if the following two conditions are satisfied:

\medskip\noindent
{\rm (i)}
$ \TC(X_{s}/k(s), p-1) $ is true for some
$ s \in |U| $.

\medskip\noindent
{\rm (ii) (0.4)} is true for
$ Y \subset X_{K}$ over
$ K $.
}

\medskip
Here we choose an embedding
$ \ok \to \oK $ so that
$ G_{k} $ is identified with a quotient of
$ G_{K} := \Gal(\oK/K) $,
because
$ K \cap \ok = k $.
We may assume
$ s \in U(k) $ replacing
$ k $ with a finite Galois extension (using the Galois action)
if necessary.
For
$ s \in U(k) $ we have the canonical isomorphism
$$
(H^{2p}(X_{\ok},\bQ_{l}(p))^{X_{s}})^{G_{k}} =
(H^{2p}(X_{\oK},\bQ_{l}(p))^{Y})^{G_{K}}.
$$
In particular, the Tate conjecture
$ \TC(X/k,p) $ implies condition~(ii) by the remark after (0.4).
For the proof of (0.7), we use the Leray spectral sequence
$$
{E}_{2}^{i,j} = H^{i}(S_{\ok}, R^{j}\of_{*}\bQ_{l})
\Rightarrow H^{i+j}(\tX_{\ok}, \bQ_{l}),
$$
which is compatible with the Galois action, and degenerates at
$ E_{2} $.
The Hodge analogue of (0.7) is given in [26], [28],
see also [22] and [41], etc.

In the case of divisors (i.e.
$ p = 1) $ with
$ n = 2 $,
S.~Bloch and K.~Kato informed us that the conjecture~(0.4) in the case
$ k $ is a number field is closely related with the finiteness of the
$ l $-primary torsion part of the Tate-Shafarevich group.
Then J.~Nekovar told us that it is also related with the de
Rham conjecture for nonproper varieties [15] after I explained him the
construction of
$ e(c) $, see also [24].
Indeed, using the Kummer sequence, we get the isomorphism
$$
H^{1}(D_{\ok}, \bZ_{l}(1)) = T_{l}J_{D}(\ok),
$$
and similarly for
$ X $, where
$ J_{D} $ denotes the Picard variety of
$ D = Y_{s} $,
and
$ T_{l} $ the Tate module.
Let
$$
J_{D,X} = \Coker(J_{X} \to J_{D}),\quad E_{D,X} =
H^{1}(D_{\ok},
\bZ_{l}(1))_{X},
$$
and
$ V_{D,X} = E_{D,X}\otimes_{\bZ_{l}}\bQ_{l} $,
so that
$$
E_{D,X} = T_{l}J_{D,X}(\ok).
$$
Then we have the short exact sequence
$$
0 \to J_{D,X}(k)\otimes \bZ_{l}\to H^{1}(G_{k},
E_{D,X}) \to T_{l}H^{1}(G_{k}, J_{D,X}(\ok)) \to 0,
\leqno(0.8)
$$
using Galois cohomology, and the conjecture~(0.4) is equivalent to
the assertion that
$ e(c) $ belongs to the image of
$ J_{D,X}(k)\otimes \bQ_{l} $ by the first morphism of (0.8).

As remarked by Bloch, Kato and Nekovar, the conjecture~(0.4) is
thus reduced to the conjecture on the finiteness of the
$ l $-primary torsion part of the Tate-Shafarevich group of
$ J_{D,X} $ and the de Rham conjecture for nonproper varieties,
using the theory of Bloch-Kato on
$ {H}_{g}^{1} $ ([3], Remark before 3.8), see (3.5) below.
(In the case $k$ is a finite field, the relation between
the Tate conjecture for $X$ and the finiteness of the Tate-Shfarevich
group of the Jacobian of the generic fiber of $X/S$ was proved by
Tate [36] where $S$ is a curve over which $X$ is dominant.)

The above two conjectures are, however, {\it still} insufficient to
deduce the Tate conjecture in our case.
Indeed, let
$$
J_{X_{s},X} = \Coker(J_{X}\otimes_{k}k(s) \to J_{X_{s}}),
\quad J_{Y,X_{K}} = \Coker(J_{X}\otimes_{k}K \to J_{Y}),
$$
and
$ \cL_{\bQ_{l}} $ the quotient of
$ R^{1}f_{*}\bQ_{l}(1)|_{U} $ by the geometrically constant part so
that
$ \cL_{\bQ_{l}} $ is identified with
$ (T_{l}J_{Y,X_{K}})\otimes_{\bZ_{l}}\bQ_{l} $, see (3.1).
We define
$ \te(c) \in \Ext^{1}(\bQ_{l}, \cL_{\bQ_{l}}) $ for
$ c \in (H^{2}(X_{\ok}, \bQ_{l}(1))^{X_{s}})^{G_{k}} $ using the
sheaf version of the short exact sequence (0.2).
(Note that
$ H^{2}(X_{\ok}, \bQ_{l}(1))^{X_{s}} $ is independent of
$ s \in |U| $.)
Then the restriction of
$ \te(c) $ to
$ s \in |U| $ coincides with
$ e(c) $ in (0.4).
Here the stalk of
$ \cL_{\bQ_{l}} $ at a geometric point over
$ s $ is
$ V_{X_{s},X} $ by the invariant cycle theorem.
Applying the above remarks of Bloch, Kato and Nekovar to the smooth closed
fibers
$ X_{s} $ of the Lefschetz pencil, and assuming the conjectures mentioned
there, we would get
$ \zeta_{s} \in J_{X_{s},X}(k(s))\otimes \bQ_{l} $ whose image by
the Abel-Jacobi map coincides with the restriction of
$ \te(c) $ to
$ s $ (in the case
$ k $ is a number field).
But it is still unclear whether the
$ \zeta_{s} $ for
$ s \in |U| $ determine an element of
$ J_{Y,X_{K}}(K)\otimes \bQ_{l} $.
Using exact sequences similar to (0.8) for
$ s \in |U| $ together with the canonical morphism of short exact
sequences, this problem is equivalent to

\medskip\noindent
{\bf 0.9.~Conjecture.}
$ T_{l}H^{1}(G, J_{Y,X_{K}}(\tK)) \to \prod_{s\in |U|}
T_{l}H^{1}(G_{s}, J_{X_{s},X}(\ok)) $ is injective.

\medskip

Here
$ G = \Gal(\tK/K) $ with
$ \tK $ the maximal subfield of
$ \oK $ that is unramified over
$ U $,
and
$ G_{s} = \Gal(\ok/k(s)) $.
(Actually the Tate conjecture is equivalent to
$ T_{l}H^{1}(G, J_{Y,X_{K}}(\tK)) = 0 $, see Remark~(3.4)(ii).)
Set
$ E = \Gamma (\Spec \tK, \cL), V = E\otimes_{\bZ_{l}}\bQ_{l} $.
As a much weaker (and easier) version of (0.9), we have at least the
injectivity
of
$$
H^{1}(G, V) \to \mprod_{s\in |U|} H^{1}(G_{s}, V).
\leqno(0.10)
$$

This is indispensable for not loosing information by taking the
restrictions to the closed points of
$ U $
(see [26] for the Hodge case).
This injectivity is informed from A.~Tamagawa in a more general case,
using Hilbert's irreducibility theorem and the theory of Frattini
subgroups (see also [30]).
We are also informed that a similar argument was used in [38].
In our case, however, this injectivity follows almost immediately from
arguments in [33], 9.1, see Remark~(3.2)(iii) below.
Furthermore we can prove (see (3.3)):

\medskip\noindent
{\bf 0.11.~Proposition.} {\it
There exists a thin subset
$ \Sigma $ of
$ U(k) $ in the sense of {\rm [33]} such that
$ H^{1}(G, E) \to H^{1}(G_{s}, E) $ is injective for
$ s \in U(k) \setminus \Sigma $.
}

\medskip
This is an analogue of N\'eron's injectivity theorem (see [21], [33]).
As a corollary of (0.11), we can solve (0.9)
using exact sequences similar to (0.8), if
$$
\rank\, J_{Y,X_{K}}(K) = \rank\, J_{X_{s},X}(k) \,\, \text{for some}\,\,
s \in U(k) \setminus \Sigma.
\leqno(0.12)
$$
Note that the last condition is not satisfied for certain
elliptic surfaces over
$ \bP^{1} $ (see [4] assuming Selmer's conjecture [29]).
However, this might occur only in the isotrivial case,
assuming some conjectures in the analytic number theory,
see Appendix of [5] for details.

Part of this work was done during my stay at the Max-Planck-Institut
f\"ur Mathematik in 1993, and I thank the staff of the institute for the
hospitality.
I also thank S.~Bloch, Y.~Ihara, K.~Kato, J.~Nekovar and A.~Tamagawa for
useful discussions, and S.~Morel and the referee for pointing out some
errors in earlier versions.

In Section 1 we explain some basic facts from the theory of
$ l $-adic mixed perverse sheaves.
Using this we prove (0.5), (0.7) in Section 2.
The divisor case is treated in Section 3.

In this paper, a variety means a separated scheme of finite type over a
field.

\newpage
\centerline{{\bf 1. Mixed Perverse Sheaves}}

\bigskip\noindent
Since the theory of {\it mixed} perverse sheaves is presented in [2] only
for
varieties over a finite field, we give a short account in the case of
varieties
$ X $ defined on a finitely generated field
$ k $ over
$ \bQ $ (see [8], [12], [19] for the case of
$ X = \Spec k $).
In this paper we restrict to the characteristic zero case.
In the case
$ k $ is finitely generated over a finite field, [2] would be essentially
sufficient
taking a formula similar to (1.16.3) for definition.

\medskip\noindent
{\bf 1.1.~Definition.} Let
$ k $ be a field of characteristic zero,
$ \ok $ an algebraic closure of
$ k $,
and
$ G_{k} = \Gal(\ok/k) $.
Let
$ X $ be a variety over
$ k $,
and
$ l $ a prime number.
Let
$ \Perv(X_{\ok}, \bQ_{l}) $ be the category of perverse sheaves on
$X_{\ok}$ with $\bQ_l$-coefficients, see [2]. 2.2.18.
Let
$ \Perv(X_{\ok},\bQ_{l};G_k) $ denote the category of objects of
$ \Perv(X_{\ok}, \bQ_{l}) $ endowed with an action of
$ G_{k} $, i.e. an object
$ \cF $ of
$ \Perv(X_{\ok},\bQ_{l};G_k) $ consists of
$ (\cF_{\ok}, u) $,
where
$ \cF_{\ok} \in \Perv(X_{\ok}, \bQ_{l}) $ and $u$ is a collection of
isomorphisms
$$
u(\gamma ) : \gamma^{*}\cF_{\ok} \to \cF_{\ok}\quad
\text{for}\,\,\, \gamma \in G_{k},
$$
satisfying the compatibility
$$
u(\gamma \gamma ') = u(\gamma )_{\circ }\gamma^{*}u(\gamma').
\leqno(1.1.1)
$$
Here
$ \gamma $ denotes also the (contravariant) action of
$ \gamma \in G_{k} $ on
$ X_{\ok}=X\otimes_k\ok $ defined by the base change.

Similarly we denote by
$ \Sh_{\sm}(X_{\ok},\bQ_{l};G_k) $ the category of (\'etale) {\it smooth}
$ \bQ_{l} $-sheaves on
$ X_{\ok} $ with an action of
$ G_{k} $ as above.
Here a smooth sheaf means that it corresponds to an $l$-adic
representation, i.e. a continuous morphism
$ \pi_{1}(X_{\ok}) \to \Aut(V)$ where $V$ is a finite dimensional
$\bQ_l$-vector spaces, see e.g. [23].
(The base point of the fundamental group is an algebraic closure
of the function field unless otherwise stated.)

For
$ \cF = (\cF_{\ok}, u) \in \Perv(X_{\ok},\bQ_{l};G_k) $,
$ \cF_{\ok} $ will be called the underlying perverse sheaf on
$ X_{\ok} $.
We have the same for
$ \cL = (\cL_{\ok}, u) \in \Sh_{\sm}(X_{\ok},\bQ_{l};G_k) $.

\medskip\noindent
{\bf 1.2.~Remarks.} (i) In the above definition, there is no condition
on the continuity of the action of $G_k$.
Later we will consider full subcategories of
$\Perv(X_{\ok},\bQ_{l};G_k)$ where the continuity of the action follows
from other conditions.
(See [2], 5.1.2 for the finite field case.)
This is considered in order to prove an analogue of [2], 5.3.8
(see (1.12)(ii) below).

\medskip
(ii) If
$ X $ is smooth and pure dimensional, let
$ \Sh_{\sm}(X,\bQ_{l}) $ denote the category of $l$-adic smooth
sheaves on $X$ (which correspond to $l$-adic representations of
$\pi_1(X)$).
Then we have fully faithful functors
$$
\Sh_{\sm}(X,\bQ_{l}) \to \Sh_{\sm}(X_{\ok},\bQ_{l};G_k)\to
\Perv(X_{\ok},\bQ_{l};G_k),
\leqno(1.2.1)
$$
where the last functor associates
$ \cL[\dim X] $ to
$ \cL \in \Sh_{\sm}(X_{\ok},\bQ_{l};G_k) $.
The full faithfulness of the first functor is proved by using
$$
H^{0}(X,\cHom(\cL,\cL'))=H^0(X_{\ok},\cHom(\rho^*\cL,\rho^*\cL'))^{G_k},
$$
where $\cL,\cL'$ are smooth sheaves and
$\rho:X_{\ok}\to X$ is the canonical morphism.
(This will be used for example in (1.12)(ii).)

\medskip
(iii) The forgetful functor
$$
\Perv(X_{\ok},\bQ_{l};G_k)\to \Perv(X_{\ok}, \bQ_{l})
\leqno(1.2.2)
$$
is exact and faithful.
This induces the forgetful functor
$$
D^{b} \Perv(X_{\ok},\bQ_{l};G_k)\to D^{b} \Perv(X_{\ok},
\bQ_{l})\to {D}_{c}^{b}(X_{\ok}, \bQ_{l}),
\leqno(1.2.3)
$$
where the last functor is given in [2], 3.1.9.

\medskip
(iv) Let
$ X_{\red} $ be the reduced variety associated with
$ X $.
Then we have an equivalence of categories
$$
\Perv(X_{\ok},\bQ_{l};G_k)= \Perv((X_{\red})_{\ok},\bQ_{l};G_k).
\leqno(1.2.4)
$$

\medskip\noindent
{\bf 1.3.~Proposition.} {\it
We have the canonical functors
$ f_{*} $,
$ f_{!} $,
$ f^{*} $,
$ f^{!} $,
$ \psi_{g} $,
$ \varphi_{g} $,
$ \bD $,
$ \boxtimes $,
$ \otimes $,
$ \cHom $ between the bounded derived categories
$ D^{b} \Perv(X_{\ok},\bQ_{l};G_k) $ for morphisms
$ f $ of
$ k $-varieties
and functions
$ g $.
Furthermore, these functors commute with the forgetful functor
{\rm (1.2.3)}.
}

\medskip\noindent
{\it Proof.}
The category of objects of
$ {D}_{c}^{b}(X_{\ok}, \bQ_{l}) $ with an action of
$ G_{k} $ as in (1.1.1) is stable by
$ \otimes $,
$ \cHom $, and the pull-back by projections.
So we get the assertion on the external product
$ {\boxtimes } $ and the dual $\bD $.

For the direct images, we can apply the same argument as in [1] if
we have the direct images of perverse sheaves with Galois action
for the embedding of the complement of a Cartier divisor, see [27], 2.4.
But the stability for this direct image easily follows from the
above definition of perverse sheaves with Galois action.

The pull-backs can be defined as the adjoint functor of the direct images,
and their existence is reduced to the case of a closed embedding where
we can use Cech-type complexes associated with an affine open covering of
the complement of the image, see [27], 3.3.

For the nearby and vanishing cycle functors, we can apply a generalization
of Deligne's construction [10] which uses finite determination sections,
see [27], 5.2 and also Remark~(1.4)(ii) below.

The last two functors are expressed using other functors.

\medskip\noindent
{\bf 1.4.~Remarks.} (i) Let
$ g $ be a function on
$ X $,
and put
$ Y = g^{-1}(0), X' = X \setminus Y $.
Let
$ \Perv((X',Y,g)_{\ok},\bQ_{l};G_k) $ be the category whose objects are
$ (\cF', \cF'', u, v) $ where
$ \cF' \in \Perv(X'_{\ok}, \bQ_{l};G_k) $,
$ \cF'' \in \Perv(Y_{\ok}, \bQ_{l};G_k) $, and
$$
u : \psi_{g,1}\cF' \to \cF'', \quad v : \cF'' \to \psi_{g,1}\cF'(1),
$$
such that
$ vu = N $.
(Here
$ \psi_{g,1}, \varphi_{g,1} $ denote the unipotent monodromy part of
$ \psi_{g}, \varphi_{g} $.)
Then, by the Deligne-MacPherson-Verdier type extension theorem [39],
we have an equivalence of categories
$$
\Perv(X_{\ok},\bQ_{l};G_k) \simto \Perv((X',Y,g)_{\ok}, \bQ_{l};G_k),
\leqno(1.4.1)
$$
induced by the functor
$$
\cF \mapsto (\cF|_{X'}, \varphi_{g,1}\cF, \can, \Var).
$$
Indeed, an inverse functor is constructed by using the functor
$ \xi_{g} $, see [25], 2.28.

Let
$ \Perv(X_{\ok},\bQ_{l};G_k)_{X'}^{\sm} $ (resp.
$ \Perv((X',Y,g)_{\ok}, \bQ_{l})_{X'}^{\sm} $) be the full subcategory of
$ \Perv(X_{\ok},\bQ_{l};G_k) $ (resp.
$ \Perv((X',Y,g)_{\ok}, \bQ_{l}) $) defined by the condition:
$ \cF|_{X'} $ (resp.
$ \cF' $) is smooth.
Then (1.4.1) induces an equivalence of categories
$$
\Perv(X_{\ok},\bQ_{l};G_k)_{X'}^{\sm}\simto\Perv((X',Y,g)_{\ok},
\bQ_{l})_{X'}^{\sm}.
\leqno(1.4.2)
$$
So an object of
$ \Perv(X_{\ok},\bQ_{l};G_k) $ can be obtained by gluing smooth sheaves
inductively.

\medskip
(ii) We define the nearby and vanishing cycle functors
$ \psi_{g}, \varphi_{g} $ so that they preserve perverse sheaves (i.e.,
they
correspond to
$ \boR\Psi [-1], \boR\Phi [-1] $ in [10]).
With the notation of Remark~(i) above, let
$ j : X' \to X $ denote the inclusion, and
$ g' : X' \to S' = \Spec k[t, t^{-1}] $ the restriction of
$ g $.
Let
$ E_{i} \,(i \ge 0) $ be a standard inductive system of indecomposable
smooth
sheaves on
$ S' $ with a weight filtration
$ W $ such that
$ {\Gr}_{j}^{W}E_{i} = \bQ_{l,S'}(-k) $ for
$ j = 2k $ with
$ 0 \le k \le i $ and
$ 0 $ otherwise.
These are constructed geometrically,
see [27], 5.1.
Then we have by definition ([27], 5.2)
$$
\psi_{g,1}\cF' = \Ker(j_{!}(\cF'\otimes g^{\prime *}E_{i}) \to
j_{*}(\cF'\otimes g^{\prime *}E_{i}))\quad\text{if}\,\,\,i\gg 0.
\leqno(1.4.3)
$$
For
$ \varphi_{g,1}\cF $ we take the cohomology of the single complex
associated to
$$
\CD
j_{!}j^{*}\cF @>>> \cF
\\
@VVV @VVV
\\
j_{!}(j^{*}\cF\otimes g^{\prime *}E_{i}) @>>>
j_{*}(j^{*}\cF\otimes g^{\prime *}E_{i}).
\endCD
$$

\medskip\noindent
{\bf 1.5.~Definition.} Assume
$ k $ is finitely generated over
$ \bQ $.
Let $\Perv(X/k,\bQ_{l})$ denote the full subcategory of
$\Perv(X_{\ok},\bQ_{l};G_k) $ consisting of objects
$ \cF $ satisfying the following condition
(by increasing induction of $\dim X$):

For each point
$ x $ of
$ X $,
there exists a Zariski-open neighborhood
$ U $ of
$ x $ in
$ X $ together with a function
$ g $ on
$ U $ such that
$ Y' := g^{-1}(0) $ has dimension
$ < \dim X $,
$ U' := U \setminus Y' $ is pure dimensional,
$ \cF|_{U'}[-\dim U'] $ is the image of $\cL\in\Sh_{\sm}(U',\bQ_l)$, and
$ \varphi_{g,1}\cF|_{U} \in \Perv(Y'/k, \bQ_{l}) $.

We define a full subcategory
$\Perv(X/k,\bQ_{l})^{\gnr} $ of
$\Perv(X/k,\bQ_{l}) $ consisting of perverse sheaves
{\it generically unramified} over $k$ by the following conditions:

In the case
$ X $ is pure dimensional and
$ \cL := \cF[- \dim X] $ is a smooth sheaf,
there exist a finitely generated
$ \bZ[\frac{1}{l}] $-subalgebra
$ R $ of
$ k $ whose fractional field is
$ k $ and an
$ R $-scheme
$ X_{R} $ of finite type whose generic fiber over
$ k $ is isomorphic to
$ X $ and such that the
$ l $-adic representation corresponding to
$ \cL $ is unramified over
$ X_{R} $, i.e. it factors through
$ \pi_{1}(X_{R}) $.

In general, there exists for each
$ x\in X $ a Zariski-open neighborhood
$ U $ of
$ x $ in
$ X $ together with a function
$ g $ on
$ U $ such that
$ Y' := g^{-1}(0) $ has dimension
$ < \dim X $,
$ U' := U \setminus Y' $ is pure dimensional,
$ \cF|_{U'}[-\dim U'] $ is a smooth sheaf and is generically
unramified over
$ k $,
and
$ \varphi_{g,1}\cF|_{U} \in \Perv(Y'/k, \bQ_{l})^{\gnr} $
(by increasing induction on $\dim X$).

\medskip\noindent
{\bf 1.6.~Remarks.} (i) Let
$ A $ a complete discrete valuation ring with the maximal ideal
$ m $ such that
$ A/m $ is a finite field of characteristic
$ l $.
Set
$ A_{i} = A/m^{i+1} $,
and let
$ K $ be the fractional field of
$ A $.
Let
$ (M_{i})_{i\in \bN} $ be a projective system of complexes of
$ A_{i} $-modules such that
$ M_{i}^{j} = 0 $ for
$ j \gg 0 $ (independently of 
$ i $),
$ H^{j}M_{0} $ are finite
$ A_{0} $-modules, and the transition morphism
$ M_{i+1} \to M_{i} $ induces an isomorphism in
$ D^{-}(A_{i}) $
$$
M_{i+1}{\otimes }_{{A}_{i+1}}^{\boL}A_{i} \simto M_{i}.
\leqno(1.6.1)
$$
Then, as in [18], pp.~474--480, there exist a complex of finite
$ A $-modules
$ L $ which is bounded above and isomorphisms
$ L_{i} := L\otimes_{A}A_{i} \simeq M_{i} $ in
$ D^{-}(A_{i}) $ which are compatible with the isomorphisms (1.6.1)
and the canonical isomorphisms
$ L_{i+1}\otimes_{A_{i+1}}A_{i} = L_{i} $.
(Here we can take
$ L $ such that the differential of
$ L_{0} $ is zero.)

Indeed, for a nonnegative integer
$ i $,
set
$ R = A_{i}, R' = A_{i+1} $ in order to simplify the notation.
Then it is enough to show the following:

Let
$ u : N \to L$, $v : N \to M $ be quasi-isomorphisms of
complexes of flat
$ R $-modules which are bounded above, and
$ M' $ a complex of flat
$ R' $-modules which is bounded above.
Assume the components of
$ L $ are finite free over
$ R $,
and we have an isomorphism
$ M'\otimes_{R'}R = M $.
Then there exist complexes of flat
$ R' $-modules
$ L', N' $ which are bounded above, together with morphisms
$ u : N' \to L'$, $v : N' \to M' $ and isomorphisms of
complexes
$ L'\otimes_{R'}R = L$, $N'\otimes_{R'}R = N $,
such that
$ u'\otimes_{R'}R $ and
$ v'\otimes_{R'}R $ are identified with
$ u $ and
$ v $ up to homotopy.

This formulation may be slightly different from loc.~cit.
However the argument is essentially the same.
Indeed, we may assume the components of
$ N, M' $ are projective over
$ R, R' $ by taking resolution, and furthermore
$ u, v $ are componentwise surjective by replacing
$ N $.
Let
$ K = \Ker\,u $.
Then
$ K $ is acyclic and
$ N $ is identified with the mapping cone of
$ M[-1] \to K $.
So the remaining argument is similar to loc.~cit.

\medskip
(ii) The theory of complexes of
$ l $-adic sheaves has become available in a general situation by
T.~Ekedahl, O.~Gabber and U.~Jannsen, see [13].
For our purpose, we can take the following formulation
(which seems more down to the earth) by modifying some
of the arguments in [7].

Let
$ X $ be a noetherian scheme on which
$ l $ is invertible, and
$ A $,
$ m $,
$ A_{i} $,
$ K $ be as in Remark~(i) above.
Let
$ M(X,A) $ denote the abelian category of projective systems
$ (\cF_{i})_{i\in \bN} $,
where
$ \cF_{i} $ are \'etale sheaves of
$ A_{i} $-modules.
Let
$ C(X,A) $ be the category of complexes of
$ M(X,A) $,
and define
$ K(X,A), D(X,A) $ using homotopy and quasi-isomorphism as in [40].
Similarly we define
$ C^{*}(X,A), K^{*}(X,A), D^{*}(X,A) $ for
$ * = +, -, b $,
so that
$ D^{*}(X,A) $ is naturally equivalent to a full subcategory of
$ D(X,A) $ (which is defined by a cohomological boundedness condition),
using
the truncations
$ \tau_{\le n}, \tau_{\ge n} $.

We say that
$ (\cF_{i}) \in $
$ D^{b}(X,A) $ is {\it strictly constructible}, if
$ \cF_{0} $ has
$ A_{0} $-constructible cohomology and the transition morphism
$ \cF_{i+1} \to \cF_{i} $ induces an isomorphism as in (1.6.1) with
$ M $ replaced by
$ \cF $ in
$ D^{b}(X,A_{i}) $, see [7], 1.1.2.
We will denote by
$ {D}_{c}^{b}(X,A) $ the full subcategory of
$ D^{b}(X,A) $ whose objects are strictly constructible.
Then
$ {D}_{c}^{b}(X,A) $ has the truncation
$ \tau '_{\le n} $ in loc.~cit. as follows:

In the usual definition of
$ \tau_{\le n}\cF_{i} $,
we replace
$ \Ker\,d \subset {\cF}_{i}^{n} $ with the subsheaf
$ {\cK}_{i}^{n} $ of
$ \Ker\,d $,
containing
$ \Imm\, d $,
such that
$$
\Imm({\cK}_{i}^{n} \to \cH^{n}\cF_{i}) = \Imm(\cH^{n}\cF_{j}
\to \cH^{n}\cF_{i})\quad \text{for }j \gg i.
$$
Here the right-hand side is independent of
$ j \gg i $ by the strict constructibility of
$ (\cF_{i}) $.
Since
$ (\tau '_{\le n}\cF_{i}) \in D(X,A) $ is well defined for
$ (\cF_{i}) \in {D}_{c}^{b}(X,A) $,
it is enough to verify (1.6.1) for
$ (\tau '_{\le n}\cF_{i}) $.
By Remark~(i) above, we may replace the stalk
$ (\cF_{i,x}) $ at each geometric point
$ x $ with
$ (L_{i}) $,
where
$ L_{i} = L\otimes_{A}A_{i} $ with
$ L $ a bounded complex of finite free
$ A $-modules.
Then
$$
H^{n}L = \varprojlim H^{n}L_{i}
$$
by the Mittag-Leffler condition, and
$$
\Imm(\varprojlim H^{n}L_{j} \to H^{n}L_{i}) =
\Imm(H^{n}L_{j} \to H^{n}L_{i})\quad \text{for }j \gg i,
$$
by the finiteness of
$ H^{n}L_{j} $ for any
$ j $.
This means
$$
\tau '_{\le n}L_{i} = (\tau_{\le n}L)\otimes_{A}A_{i}.
$$
Thus
$ (\tau '_{\le n}\cF_{i}) $ is strictly constructible.

Let
$ \cC $ denote the heart of this
$ t $-structure.
Then
$ \cC $ is naturally equivalent to
$ \Sh_{c}(X,A) $ the category of
$ A $-constructible sheaves on
$ X $ by [7], 1.1.2.
Indeed, we have naturally
$ \alpha : \cC \to \Sh_{c}(X,A) $ by
$ \alpha (\cF_{i}) = (\cH^{0}\cF_{i}) $.
For
$ \beta : \Sh_{c}(X,A) \to \cC $,
let
$ (\cS_{i}) \in \Sh_{c}(X,A) $ with
$ \cS_{i} = \cS_{j}\otimes_{A_{j}}A_{i} $ for
$ j > i $,
and take a quasi-isomorphism
$ (\cE_{i}) \to (\cS_{i}) $ in
$ C^{-}(X,A) $ such that the stalks of the components of
$ \cE_{i} $ at geometric points are flat over
$ A_{i} $.
Let
$ k $ be a positive integer such that the torsion of
$ \varprojlim \cS_{i,x} $ at geometric points
$ x $ is annihilated by
$ m^{k} $.
Then we have as in loc.~cit.
$$
\beta (\cS_{i}) = (\tau_{\ge -1}(\cE_{i+k}\otimes_{A_{i+k}}A_{i})).
$$
Indeed,
$ \alpha \beta \simeq id $ is clear, and
$ \beta \alpha \simeq id $ is induced by
$ (\cF_{i}) \to (\cH^{0}\cF_{i}) $ for
$ (\cF_{i} ) \in \cC $,
combined with the quasi-isomorphisms
$ \cF_{i+k}\otimes_{A_{i+k}}A_{i} \to \cF_{i} $ where we may assume
that the stalks of the components of
$ \cF_{i} $ are flat over
$ A_{i} $.

\medskip
(iii) Let
$ D_{c}^{b}(X_{R}, \bQ_{l}) $ be the derived category of bounded
complexes of
$ l $-adic sheaves on
$ X_{R} $ with constructible cohomology (see Remark~(ii) above),
where
$ R $ and
$ X_{R} $ are as in (1.5).
We have a canonical functor
$$
\iota_{R,\ok} : D_{c}^{b}(X_{R}, \bQ_{l}) \to
D_{c}^{b}(X_{\ok}, \bQ_{l}; G_{k}),
$$
where the target is the category consisting of objects of
$ D_{c}^{b}(X_{\ok}, \bQ_{l}) $ with an action of
$ G_{k} $ as in (1.1).
(This category is considered only to define
$ \Perv(X/k,\bQ_{l})^{\gen} $ in Remark~(v) below, and
it is not clear whether an object of
$ D_{c}^{b}(X_{\ok}, \bQ_{l}; G_{k}) $ is represented by
a complex with an action of
$ G_{k} $ in
$ C(X_{\ok}, \bQ_{l}) $.)
The functor
$ \iota_{R,\ok} $ commutes with direct images and pull-backs
by the generic base change theorem [11], 2.9.

\medskip
(iv) The category $\Perv(X,\bQ_l)$ of $l$-adic perverse sheaves on
$X$ can be defined as a full subcategory of $D_c^b(X,\bQ_l)$ using [2]
since the latter category is defined as in Remark~(ii) above.
We have the gluing of perverse sheaves as in (1.4.1) by the same
argument, and hence this category is equivalent to
$ \Perv(X/k,\bQ_l)$ by increasing induction on $\dim X$.
(Note that perverse sheaves can be defined locally, see [2], 3.2.4).
As a corollary, (1.3) holds also for $D^b\Perv(X/k,\bQ_l)$.
If the reader prefers, the category $\Perv(X/k,\bQ_l)$ in (1.5) may
be defined by the above $\Perv(X,\bQ_l)$
(however, the inductive argument using (1.4.1) will be needed to
define certain full subcategories of $\Perv(X/k,\bQ_l)$).

By Remark~(iii) above, a similar argument applies to
$ \Perv(X/k,\bQ_{l})^{\gnr} $, and this category is stable by the
cohomological functors associated with the functors in (1.3), i.e.,
$ H^{i}f_{*} $, etc.
Here the stability by dual and external product follows from the
commutativity of these functors with
$ \varphi_{g,1} $.
Then the assertion of (1.3) holds also for
$ D^b\Perv(X/k,\bQ_{l})^{\gnr} $,
and the functors in (1.3) commute with the forgetful functors
$$
D^{b} \Perv(X/k,\bQ_{l})^{\gnr} \to
D^{b} \Perv(X/k,\bQ_{l}) \to
D^{b}\Perv(X_{\ok},\bQ_{l};G_k).
$$

We have the stability of
$ \Perv(X/k,\bQ_{l})^{\gnr} $ and
$ \Perv(X/k,\bQ_{l}) $ by
subquotients in
$ \Perv(X_{\ok},\bQ_{l};G_k) $, because
$ \varphi_{g,1} $ is an exact functor.
Indeed, in case of smooth sheaves, a subobject of the image of a
smooth sheaf $\cL$ in $\Sh_{\sm}(X_{\ok},\bQ_l;G_k)$ defines a subsheaf
of $\cL$ by reducing to the finite coefficient case
(because any $l$-adic representation of a profinite group $G$ on a
$\bQ_l$-vector space $V$ has a finitely generated $\bZ_l$-submodule
which generates $V$ over $\bQ_l$ and is stable by the action of $G$)
and using the action of $G_k$ on the variety representing the torsion
sheaf over $X_{\ok}$
(because the stability of a subvariety by the Galois action means
that the subvariety is defined over $k$).
This will be used for example in (1.12)(ii).

\medskip
(v) Let
$ \Perv(X/k,\bQ_{l})^{\gen} $ denote the full subcategory of
$ \Perv(X_{\ok},\bQ_{l};G_k) $ consisting of objects whose image in
$ D_{c}^{b}(X_{\ok}, \bQ_{l}; G_{k}) $ belongs to the
essential image of the functor
$ \iota_{R,\ok} $ in Remark~(iii) above for some
$ R $ (which may depend on the objects).
Then
$ \Perv(X/k,\bQ_{l})^{\gen} $ is stable by cohomological direct
images and pull-backs by [11], 2.9.
(Here we may assume that every stratum of the stratification of
$ X_{R} $ associated with a complex is dominant over
$ \Spec R $ replacing
$ R $ if necessary.)
Hence it is also stable by vanishing cycle functors.
By induction on
$ \dim X $ and using the gluing as in (1.4.1),
we get an equivalence of categories
$$
\Perv(X/k,\bQ_{l})^{\gen} = \Perv(X/k,\bQ_{l})^{\gnr}.
$$

\medskip\noindent
{\bf 1.7.~Definition.} Let
$ Z $ be an irreducible closed subvariety of
$ X $.
We say that
$ \cF \in \Perv(X/k,\bQ_{l}) $ has {\it strict support}
$ Z $,
if
$ \supp \cF = Z $ or
$ \emptyset $,
and if
$ \cF $ has no nontrivial sub or quotient objects with smaller support.
We will denote by
$ \Perv(X/k,\bQ_{l})_{Z} $ the full subcategory of
$ \Perv(X/k,\bQ_{l}) $ consisting of objects with strict support
$ Z $.

We say that
$ \cF \in \Perv(X/k,\bQ_{l}) $ admits {\it a strict support
decomposition}, if we have a decomposition
$$
\cF = \moplus_{Z} \cF_{Z}\quad \text{with\quad } \cF_{Z} \in
\Perv(X/k, \bQ_{l})_{Z}.
\leqno(1.7.1)
$$
We will call
$ \cF_{Z} $ the direct factor of
$ \cF $ with strict support
$ Z $.
We have the same for
$ \Perv(X_{\ok}, \bQ_{l}) $.

\medskip\noindent
{\bf 1.8.~Remarks.} (i) The strict support decomposition (1.7.1) is
unique, because
$$
\Hom(\cF_{Z}, \cF_{Z'}) = 0\quad \text{if\quad } Z \ne Z'.
\leqno(1.8.1)
$$

\medskip
(ii) Assume that
$ \cF \in \Perv(X/k,\bQ_{l}) $ admits strict support decomposition
and
$ \cF_{\ok} $ is semisimple.
Then any filtration of
$ \cF $ is compatible with the strict support decomposition (1.7.1).
Indeed, for an exact sequence
$$
0 \to \cF' \to \cF \to \cF'' \to
0
$$
in
$ \Perv(X/k,\bQ_{l}) $,
$ \cF' $ and
$ \cF'' $ admit strict support decomposition, and we have an exact sequence

$$
0 \to \cF'_{Z} \to \cF_{Z} \to \cF''_{Z}
\to 0.
$$
This can be proved using (1.8.1) (or (1.8)(iv) below).

\medskip
(iii) For a locally closed embedding
$ j : X \to Y $,
we have the intermediate direct image
$ j_{!*} $ (see [2]) defined by
$$
j_{!*}\cF = \Imm(H^{0}j_{!}\cF \to H^{0}j_{*}\cF)\quad
\text{for}\,\,\,
\cF \in \Perv(X/k,\bQ_{l}).
$$
Let
$ \cF \in \Perv(X/k,\bQ_{l})_{Z} $,
and
$ Z' $ be a dense open subvariety of
$ Z $ with the inclusion morphism
$ j : Z' \to X $.
Then we have a canonical isomorphism
$$
\cF = j_{!*}j^{*}\cF,
\leqno(1.8.2)
$$
using the morphisms
$ H^{0}j_{!}j^{*}\cF \to \cF \to H^{0}j_{*}j^{*}\cF $.
Furthermore, (1.8.2) holds with
$ \cF $ replaced by
$ \cF_{\ok} $,
and
$ j $ by its base extension over
$ \ok $.
So, if
$ \cF \in \Perv(X/k,\bQ_{l}) $ admits strict support decomposition,
then
$ \cF_{\ok} $ does.

\medskip
(iv) For
$ \cF \in \Perv(X/k,\bQ_{l}) $,
it admits strict support decomposition if and only if the composition
$$
i_{*}H^{0}i^{!}\cF \to \cF \to i_{*}H^{0}i^{*}\cF
\leqno(1.8.3)
$$
is an isomorphism for any closed embedding
$ i : Y \to X $,
because
$ i_{*}H^{0}i^{!}\cF $ (resp.
$ i_{*}H^{0}i^{*}\cF) $ is identified with the largest sub (resp. quotient)
object
of
$ \cF $ supported in the image of
$ Y $.
This is the same for
$ \Perv(X_{\ok}, \bQ_{l}) $.
Then, combining this with Remark~(iii) above,
$ \cF \in \Perv(X/k,\bQ_{l}) $ admits strict support decomposition,
if and only if
$ \cF_{\ok} $ does.
(Here we can use also (1.8.1).)

\medskip\noindent
{\bf 1.9.~Definition.} We denote by
$ \bQ_{l} $ the constant object of
$ \Perv(\Spec k/k, \bQ_{l}) $ with trivial Galois action.
For a variety
$ X $ with a canonical morphism
$ a_{X} : X \to \Spec k $, we define
$$
\aligned
\bQ_{l,X}
&= (a_{X})^{*}\bQ_{l}\in D^{b} \Perv(X/k,\bQ_{l}).
\\
H^{j}(X/k,\bQ_{l})
&=H^{j}(a_{X})_{*}\bQ_{l,X} \in \Perv(\Spec k/k, \bQ_{l}).
\endaligned
$$
More generally, we set for
$ \cF \in D^{b} \Perv(X/k,\bQ_{l}) $:
$$
H^{j}(X/k,\cF) = H^{j}(a_{X})_{*}\cF,\quad {H}_{c}^{j}(X/k,\cF) =
H^{j}(a_{X})_{!}\cF
\leqno(1.9.1)
$$
If
$ X $ is smooth and pure dimensional, then
$ \bQ_{l,X} \in \Sh_{\sm}(X,\bQ_{l}) $ and
$ \bQ_{l,X}[\dim X] \in \Perv(X/k,\bQ_{l}) $.

Let
$ X $ be a pure dimensional variety, and
$ U $ a dense smooth open subvariety of
$ X $ with the inclusion
$ j : U \to X $.
Then, using (1.2.1), we define for
$ \cL \in \Sh_{\sm}(U, \bQ_{l}) $ the intersection complex
with coefficients in
$ \cL $ by
$$
\IC_{X}\cL = j_{!*}\cL[\dim X]\in \Perv(X/k,\bQ_{l}).
$$
If
$ \cL $ is the constant sheaf
$ \bQ_{l,X} $,
we will denote
$$
\IC_{X}\bQ_{l}=j_{!*}\bQ_{l,U}[\dim X],
$$

\medskip\noindent
{\bf 1.10~Remark.} For
$ \cF \in \Perv(X/k,\bQ_{l})_{Z} $,
there exist a dense smooth open affine subvariety
$ Z' $ of
$ Z $ and
$ \cL \in \Sh_{\sm}(Z', \bQ_{l}) $ such that we have an isomorphism
$$
\cF = \IC_{Z}\cL.
\leqno(1.10.1)
$$
This follows from Remark~(1.8)(iii).
Note that, if (1.10.1) holds for
$ Z' $ and
$ \cL $,
it holds also for any dense open subvariety of
$ Z' $ and the restriction of
$ \cL $.

\medskip\noindent
{\bf 1.11.~Definition.} With the notation of (1.1), assume
$ k $ is finitely generated over
$ \bQ $.
In the case
$ X $ is smooth, we say that
$ \cL \in \Sh_{\sm}(X,\bQ_{l}) $ is generically unramified and
{\it pure} of weight
$ n $,
if there exist a finitely generated
$ \bZ[\frac{1}{l}] $-subalgebra
$ R $ of
$ k $ whose fractional field is
$ k $,
a smooth scheme
$ X_{R} $ of finite type over
$ \Spec R $ whose generic fiber is isomorphic to
$ X $,
and a smooth
$ \bQ_{l} $-sheaf
$ \cL_{R} $ on
$ X_{R} $ whose restriction to
$ X $ is isomorphic to
$ \cL $ and whose stalk at each closed point is pure of weight
$ n $ in the sense of [7], see Remark~(1.12)(i) below.

We say that
$ \cF \in \Perv(X/k,\bQ_{l})_{Z} $ is {\it pure} of weight
$ n $, if
$ \cL $ in (1.10.1) is generically unramified and pure of weight
$ n - \dim Z $.
We say that
$ \cF \in \Perv(X/k,\bQ_{l}) $ is {\it pure} of weight
$ n $ if
$ \cF $ admits the strict support decomposition (1.7.1) and each
$ \cF_{Z} $ is pure of weight
$ n $.
We say that
$ \cF \in \Perv(X/k,\bQ_{l}) $ is {\it mixed}, if
$ \cF \in \Perv(X/k,\bQ_{l})^{\gnr} $ and
$ \cF $ has a finite increasing filtration
$ W $, which is called the {\it weight filtration}, such that
each
$ {\Gr}_{n}^{W}\cF $ is pure of weight
$ n $.

We will denote by
$ \Perv(X/k,\bQ_{l})^{m} $ the full subcategory of
$ \Perv(X/k,\bQ_{l})^{\gnr} $ consisting of mixed perverse sheaves.

\medskip\noindent
{\bf 1.12.~Remarks.} (i) The condition in the pure case means that the
$ l $-adic representation of
$ \pi_{1}(X) $ corresponding to
$ \cL $ is unramified over
$ X_{R} $ (i.e., it factors through
$ \pi_{1}(X_{R})) $,
and the eigenvalues of Frobenius at each closed point
$ x $ of
$ X_{R} $ are algebraic numbers whose any conjugates over
$ \bQ $ have absolute value
$ q^{n/2} $,
where
$ q = |\kappa (x)| $.
In particular, pure perverse sheaves are stable by subquotients in
$ \Perv(X/k,\bQ_{l}) $ using Remark~(1.8)(ii).
Note that pure perverse sheaves are generically unramified over
$ k $,
and hence are mixed, since
$ \Perv(X/k,\bQ_{l})^{\gnr} $ is stable by intermediate direct images,
see Remark~(1.6)(iv).

\medskip
(ii) If
$ \cF \in \Perv(X/k,\bQ_{l})^{\gnr} $ is pure, then
$ \cF_{\ok} $ is semisimple.
The argument is essentially the same as in [2], 5.3.8.
Indeed, for a short exact sequence of pure objects of
$ \Perv(X/k,\bQ_{l})^{\gnr} $
$$
0 \to \cF' \to \cF \to \cF'' \to 0,
\leqno(1.12.1)
$$
we have a commutative diagram
$$
\CD
\Ext^{1}(\cF'', \cF') @>>> \Ext^{1}(\cF''_{\ok}, \cF'_{\ok})
\\
@VVV @VVV
\\
\Hom(\bQ_{l}, H^{1}(X/k, \cHom(\cF'', \cF'))) @>>>
\Hom(\bQ_{l}, H^{1}(X_{\ok}, \cHom(\cF''_{\ok},
\cF'_{\ok}))),
\endCD
$$
using the adjunction for
$ X \to \Spec k $.
Here the right vertical morphism is an isomorphism.
Since
$ H^{1}(X/k, \cHom(\cF'', \cF')) $ has weights
$ > 0 $ by the generic base change theorem ([11], 2.9, 2.10),
the horizontal morphisms are zero, and the exact sequence of
the underlying perverse sheaves of (1.12.1) on
$ X_{\ok} $ splits.
Then it is enough to apply this to the largest semisimple subobject of
$ \cF_{\ok} $,
which is stable by the Galois action, and determines a subobject of
$ \cF $, see (1.2)(ii) and (1.6)(iv).

Note that in the above argument, we have used only the following
condition:
\begin{itemize}
\item[$(C)$]
The mod $p$ reduction of $\cF$ over any closed point of a
sufficiently small open subvariety of $\Spec R$ in (1.11) is pure.
\end{itemize}
Here the mod $p$ reduction exists by (1.6)(v).
Condition~(C) implies thus the strict support decomposition
(1.7.1) for $\cF_{\ok}$, and hence for $\cF$, see (1.8)(iv).
So we may take this as the definition of pure perverse sheaf.

\medskip
(iii)
$ \Perv(X/k,\bQ_{l})^{m} $ is stable by subquotients in
$ \Perv(X/k,\bQ_{l}) $,
and the weight filtration on a mixed perverse sheaf
$ \cF $ is unique.
Indeed, any filtration on
$ {\Gr}_{n}^{W}\cF $ is compatible with the strict support decomposition
(1.7.1) by
Remark~(ii) above and Remark~(1.8)(ii), and induces a filtration on
each
$ ({\Gr}_{n}^{W}\cF)_{Z} $.
This implies the uniqueness of
$ W $ and the stability of mixed perverse sheaves by subquotients, see
Remark~(1.6)(iv).

\medskip
(iv) By the next proposition, the equivalence of categories (1.4.1) holds
also
for mixed perverse sheaves, because they are stable by the functor
$ \zeta_{g} $ which is used to construct the inverse functor of (1.4.1).
As a corollary, we can also define the mixed perverse sheaves as in (1.5)
by induction
on
$ \dim X $ using the functor
$ \varphi_{g} $.

\medskip\noindent
{\bf 1.13.~Proposition.} {\it
The assertion of {\rm (1.3)} holds with
$ \Perv(X_{\ok},\bQ_{l};G_k) $ replaced by
$ \Perv(X/k,\bQ_{l}) $,
$ \Perv(X/k,\bQ_{l})^{\gnr} $ and
$ \Perv(X/k,\bQ_{l})^{m} $.
}

\medskip\noindent
{\it Proof.}
The assertion for $ \Perv(X/k,\bQ_{l}) $, $ \Perv(X/k,\bQ_{l})^{\gnr} $
is shown in (1.6)(iv), and it is enough to consider the case of
$ \Perv(X/k,\bQ_{l})^{m} $.
We have the stability by the dual and external product,
because they are exact functors and the stability of pure objects
by these functors follows easily from the definition.
By the same argument as in the proof of (1.3),
it is enough to show the stability by
$ j_{*} $ when
$ j $ is the open embedding of the complement of a divisor defined by a
function
$ g $.

In this case, we have to show the existence of the relative
monodromy filtration
$ W $ on
$ \psi_{g,1}\cF $ (see [7], 1.6.13),
because $W$ defines the weight filtration on
$ j_{*}\cF $ by a formula of Steenbrink-Zucker ([34], 4.11), see
[25], 2.11.
Note that
$ j_{*}\cF $ corresponds to
$ (\cF, \psi_{g,1}\cF, N, id) $ by (1.4.1).
The existence of
$ W $ is reduced to the finite field case by restricting to the fiber over
a sufficiently general closed point of
$ \Spec R $,
and using a criterion for the existence of $W$ in [34], 2.20
(see also [25], 1.2).
Then the assertion is reduced to Gabber's result on the coincidence
of the monodromy and weight filtrations (up to a shift) in the case
$ \cF $ is pure.
(The last result does not seem to have been published, but a (possibly)
similar argument can be found in [27], 6.11.)

\medskip\noindent
{\bf 1.14.~Remark.} If
$ \cF \in \Perv(X/k,\bQ_{l})^{m} $ and
$ \cF' \in \Perv(Y/k, \bQ_{l})^{m} $ have weights
$ \ge n $ (resp.
$ \le n) $,
then
$ H^{i}f_{*}\cF, H^{i}f^{!}\cF' $ (resp.
$ H^{i}f_{!}\cF, H^{i}f^{*}\cF') $ have weights
$ \ge n + i $ (resp.
$ \le n + i) $.
For the direct images, it is enough to show the assertion for
$ H^{i}f_{*}\cF $ (using
$ f_{*}\bD = \bD f_{!}) $.
If $f$ is proper and $\cF$ is pure of weight $n$, then
$H^if_*\cF$ is pure of weight $i+n$ using the stability of pure
complexes by the direct image under a proper morphism in [7].
So the assertion is reduced to the case
$ f $ is an affine open embedding
$ j $ as above, and moreover
$ \cF $ is pure.
Then the assertion follows from the construction of
$ W $ using the monodromy filtration, see [25], 2.11.

For pull-backs, we may assume
$ X $ is a closed subvariety of
$ Y $ defined by a function (by factorizing
$ f $) since the assertion is local and the case of a smooth morphism
is easy.
Then we may assume that
$ \cF' $ is pure of weight
$ n $ with strict support and
$ \supp\cF' $ is not contained in
$ X $.
In this case we have to show that the cokernel of
$ \cF' \to j_{*}j^{*}\cF' $ has weights
$ > n $,
where
$ j $ denotes the inclusion of
$ Y \setminus X $.
But this also follows from the construction of
$ W $ using the monodromy filtration, see loc.~cit.

\medskip\noindent
{\bf 1.15.~Definition.} Let
$ S $ be an integral affine variety over
$ k $,
and
$ K = k(S) $.
For a variety
$ X $ over
$ K $,
let
$ X_{S} $ be a
$ k $-variety over
$ S $ whose generic fiber is isomorphic to
$ X $ (restricting
$ S $ if necessary).
Let
$$
\Perv(X/K/k, \bQ_{l})=\varinjlim \Perv(X_U/k, \bQ_{l})
\leqno(1.15.1)
$$
where the inductive limit is taken over nonempty open subvarieties
$ U $ of
$ S $,
and
$ X_{U} = X_{S}|_{U} $.
Similarly,
$ \Perv(X/K/k, \bQ_{l})^{\gnr} $ and
$ \Perv(X/K/k, \bQ_{l})^{m} $ are defined by replacing
$ \Perv(X_U/k, \bQ_{l}) $ with
$ \Perv(X_{U}/k, \bQ_{l})^{\gnr} $ and
$ \Perv(X_{U}/k, \bQ_{l})^{m} $ in (1.15.1).

\medskip\noindent
{\bf 1.16.~Proposition.} {\it
We have a fully faithful functor
$$
\Perv(X/K/k, \bQ_{l})\to \Perv(X/K, \bQ_{l})
\leqno(1.16.1)
$$
whose essential image is stable by subquotients in
$ \Perv(X/K, \bQ_{l}) $.
Furthermore, {\rm (1.16.1)} induces equivalences of categories
}
$$
\Perv(X/K/k, \bQ_{l})^{\gnr}
\to \Perv(X/K, \bQ_{l})^{\gnr},
\leqno(1.16.2)
$$
$$
\Perv(X/K/k, \bQ_{l})^{m}
\to \Perv(X/K, \bQ_{l})^{m}.
\leqno(1.16.3)
$$

\medskip\noindent
{\it Proof.}
To define the functor (1.16.1), we choose an embedding
$ \ok \to \oK $ which induces a morphism
$ R'\otimes_{k}\ok \to \oK $,
where
$ S = \Spec R' $.
Then (1.16.1) is given by the base change by this morphism.
Since the assertions are local, we may assume
$ X $ affine.
Let
$ g $ be a function on
$ X $,
and put
$ Y = g^{-1}(0), X' = X \setminus Y $.
We may assume
$ g $ (and
$ Y, X') $ defined over
$ S $,
shrinking
$ S $ if necessary.
Applying (1.4.1) to
$ X, X', Y $ over
$ K $ and also to
$ X_{U}, X'_{U}, Y_{U} $ over
$ k $,
and using the generic base change theorem ([11], 2.9),
the assertion is reduced to the case of smooth sheaves on smooth
varieties by induction on
$ \dim X $.
So the assertion follows.

\bigskip\bigskip\centerline{{\bf 2. Tate Conjecture}}

\bigskip\noindent
In this section,
$ k $ is a finitely generated field over
$ \bQ $, and we prove (0.5) and (0.7).
Although the arguments are similar to those in [28] in some places,
there are certain differences between the Hodge and
$ l $-adic settings, and we repeat some of the arguments here.

\medskip\noindent
{\bf 2.1.~Cycle map.} Let
$ X $ be a smooth variety over
$ k $.
With the notation of (1.9), we have a cycle map
$$
\CH^{p}(X) \to \Ext^{2p}(\bQ_{l,X}, \bQ_{l,X}(p)),
\leqno(2.1.1)
$$
where
$ \Ext $ is taken in
$ D^{b} \Perv(X/k,\bQ_{l}) $.
For a closed irreducible reduced subvariety
$ Z $ of
$ X $,
the image of the cycle
$ [Z] $ by (2.1.1) is given by the composition of the canonical
morphism
$$
\bQ_{l,X} \to \bQ_{l,Z} \to
\IC_{Z}\bQ_{l}[- \dim Z]
$$
with its dual using
$$
\bD(\IC_{Z}\bQ_{l})=\IC_{Z}\bQ_{l}(\dim Z),\quad
\bD(\bQ_{l,X}) = \bQ_{l,X}(\dim X)[2 \dim X].
$$
We have the adjunction isomorphism
$$
\Ext^{2p}(\bQ_{l,X}, \bQ_{l,X}(p)) = \Ext^{2p}(\bQ_{l},
(a_{X})_{*}\bQ_{l,X}(p)).
$$
If
$ X $ is smooth projective, we have a (noncanonical) decomposition
by [6]:
$$
(a_{X})_{*}\bQ_{l,X} \simeq \oplus H^{j}(X/k, \bQ_{l})[-j].
\leqno(2.1.2)
$$
(This holds also in the case
$ X $ is smooth proper, see [27], 6.8.)
So (2.1.1) induces naturally the cycle map
$$
cl : \CH^{p}(X) \to \Hom(\bQ_{l}, H^{2p}(X/k, \bQ_{l})(p)),
$$
which coincides with (0.1).
Let
$$
\CH_{\homm}^{p}(X) = \Ker\,cl.
$$
Then (2.1.1) induces naturally
$$
\CH_{\homm}^{p}(X) \to \Ext^{1}(\bQ_{l},H^{2p-1}(X/k,\bQ_{l})(p)),
\leqno(2.1.3)
$$
which is called the Abel-Jacobi map.
Note that the cycle map (2.1.1) can also be defined by using the Gysin
morphism

$$
(a_{Z})_{*}\bD(\bQ_{l,Z}) \to (a_{X})_{*}\bD(\bQ_{l,X})
$$
whose mapping cone is isomorphic to
$ (a_{U})_{*}\bQ_{l,X}(\dim X)[2 \dim X] $,
where
$ Z $ is the support of a cycle and
$ U = X \backslash Z $.
So (2.1.3) can be obtained also by taking the pull-back of the exact
sequence
$$
0 \to H^{2p-1}(X/k, \bQ_{l})\to H^{2p-1}(U/k, \bQ_{l}
\to {H}_{Z}^{2p}(X/k, \bQ_{l}),
$$
see [19] (and also [14], [16], [26], [28] for the Hodge case.)

\medskip\noindent
{\bf 2.2.~Remarks.} (i) Let
$ \TC(X/k,p) $ denote the Tate conjecture as in Introduction.
Then
$ \TC(X/k,p) $ depends only on
$ X $ and
$ p $,
and is essentially independent of
$ k $.
Indeed, for a finite extension
$ k' $ of
$ k $ with degree
$ d $,
$ \Spec k'\otimes_{k}\ok $ consists of
$ d $ points on which
$ \Gal(\ok/k) $ acts transitively, and its stabilizer can be
identified
with
$ \Gal(\ok/k') $.

\medskip
(ii) Let
$ k' $ be a finite extension of
$ k $.
Then
$ \TC(X/k,p) $ follows from
$ \TC(X\otimes_{k}k'/k',p) $.
This can be proved using the action of
$ \Gal(k''/k) $ on the cycles, where
$ k'' $ is a Galois extension of
$ k $ containing
$ k' $.
In particular, we may assume
$ X $ absolutely irreducible as in [35].

\medskip
(iii) More generally, let
$ k \to K $ be an arbitrary extension, and
$ \oK $ an algebraic closure of
$ K $ with an inclusion
$ \ok \to \oK $.
We have a canonical isomorphism
$$
H^{2p}(X_{\ok}, \bQ_{l})(p) = H^{2p}(X_{\oK}, \bQ_{l})(p)
\leqno(2.1.4)
$$
by the base change theorem.
Let
$ c \in H^{2p}(X_{\ok}, \bQ_{l})(p) $ which is invariant by the
Galois action.
Then, to show that
$ c $ is algebraic, it is enough to construct a cycle on
$ X_{K} $ whose cycle class coincides with
$ c $ by the isomorphism (2.1.4).
Indeed, we may assume
$ K $ finitely generated over
$ k $.
Let
$ R $ be a finitely generated algebra over
$ k $ such that
$ K $ is the fractional field of
$ R $ and the cycle is defined over
$ R $.
Then the assertion is reduced to the above Remark~(ii) by restricting to
the
fiber over a general closed point of
$ \Spec R $.

\medskip
(iv) The Tate conjecture over a
$ p $-adic field is not true.
Indeed, there exist elliptic curves without complex multiplication, but
having a formal complex multiplication, as remarked by J.-P.~Serre
([31], Remark~(i) in A.2.2).
According to Y.~Ihara, this can be verified by using a theory of Honda
[17] on formal
groups and a theory of Serre and Tate [32] on
$ p $-divisible groups.

\medskip\noindent
{\bf 2.3.}
Let
$ f : X \to S $ be a projective morphism of
$ k $-varieties.
Let
$ \cF $ be a pure object of
$ \Perv(X/k,\bQ_{l})^{m} $.
With the notation of (1.9.1), we have the Leray spectral sequence in
$ \Perv(\Spec k/k, \bQ_{l}) $:
$$
{E}_{2}^{p,q} = H^{p}(S/k, H^{q}f_{*}\cF) \Rightarrow H^{p+q}(X/k, \cF),
\leqno(2.3.1)
$$
which degenerates at
$ E_{2} $.
Indeed, the spectral sequence is induced by the filtration
$ \tau $ on
$ f_{*}\cF $,
and the
$ E_{2} $-degeneration follows from [6].
(If
$ f $ is assumed only proper, we can use [27], 6.8.)
We denote by
$ L $ the associated filtration on
$ H^{j}(X/k, \cF) $ so that
$$
{\Gr}_{L}^{i}H^{j}(X/k, \cF) = H^{i}(S/k, H^{j-i}f_{*}\cF).
\leqno(2.3.2)
$$

If
$ X $ is smooth and purely
$ n $-dimensional, and
$ \cF $ is the constant perverse sheaf
$ \bQ_{l,X}[n] $,
then we have by definition (see (1.9))
$$
H^{j}(X/k, \bQ_{l,X}[n]) = H^{j+n}(X/k, \bQ_{l}),
\leqno(2.3.3)
$$
and we denote also by
$ L $ the induced filtration on
$ H^{j+n}(X/k, \bQ_{l}) $.
Note that
$ L $ induces also the filtration
$ L $ on
$$
\Hom(\bQ_{l}, H^{2p}(X/k, \bQ_{l})(p)),
$$
because
$ L $ on
$ H^{j}(X/k, \bQ_{l}) $ splits in
$ \Perv(\Spec k/k, \bQ_{l}) $ using the decomposition theorem for the
direct image
$ f_{*}\cF $ as above.

\medskip\noindent
{\bf 2.4.}
With the above notation, assume
$ f $ is a Lefschetz pencil
$ \tX \to S $ as in Introduction, and
$ \cF = \bQ_{l,\tX}[n] $ as in (2.3.3).
Note that (0.6) is equivalent to
$$
H^{0}f_{*}(\bQ_{l,\tX}[n])_{\ok}[-1]\,\,\, \text{is an }l
\text{-adic sheaf on }S_{\ok},
\leqno(2.4.1)
$$
i.e.,
$ H^{0}f_{*}(\bQ_{l,\tX}[n]) $ has no direct factor whose
support is zero-dimensional.
Indeed, (0.6) is equivalent to the surjectivity (or nonvanishing) of
$$
\can : \psi_{t,1}H^{0}f_{*}(\bQ_{l,\tX}[n])_{\ok} \to
\varphi_{t,1}H^{0}f_{*}(\bQ_{l,\tX}[n])_{\ok}\, (= \bQ_{l}(-p)),
$$
(where
$ t $ is a local coordinate at a critical value of
$ f $), since
$ H^{n-1}(Y_{\oK}, \bQ_{l})_{X} $ is generated by the vanishing
cycles, which are locally identified with the image of the dual of
the above morphism can.
(Note that the vanishing cycles are conjugate to each other
by the global monodromy group since the discriminant in the dual
projective space parametrizing the singular hyperplane sections is
the image of a projective bundle over
$ X $, and hence is irreducible.)
So the equivalence follows from (1.4.1).

By the Picard-Lefschetz formula [9], we have:
$$
H^{j}f_{*}(\bQ_{l,\tX}[n])_{\ok}[-1]\,\,\, \text{is a
constant sheaf on }
S_{\ok}\,\,\, \text{for }j \ne 0.
$$
So we get in the notation of (2.3)
$$
{E}_{2}^{p,q} = 0\quad \text{for}\,\,\, p = 0, \,q \ne 0,
\leqno(2.4.2)
$$
because
$ S = \bP^{1} $.
Let
$ X_{s} $ be the fiber of a
$ k $-rational point
$ s $ of
$ U $ with the inclusion
$ \ti : X_{s} \to \tX $,
where
$ U $ is as in Introduction.
Then we have
$$
L^{0}H^{j}(\tX/k, \bQ_{l})
=\Ker(\ti^{*} : H^{j}(\tX/k, \bQ_{l})\to
H^{j}(X_{s}/k, \bQ_{l})),
\leqno(2.4.3)
$$
$$
L^{1}H^{j}(\tX/k, \bQ_{l})
=\Imm(\ti_{*} : H^{j-2}(X_{s}/k, \bQ_{l})(-1) \to
H^{j}(\tX/k, \bQ_{l})),
\leqno(2.4.4)
$$
where
$ \ti^{*} $ and
$ \ti_{*} $ are the restriction and Gysin morphisms.

\medskip\noindent
{\bf 2.5.~Proof of (0.5).}
By Remark~(2.2)(ii), we may assume
$ k(s) = k $ for the assertion (i) replacing
$ X $ with
$ X\otimes_{k}k(s) $ if necessary.
Then the assertion (i) follows from the weak Lefschetz theorem.
For (ii), we use the Leray spectral sequence as above.
Take
$$
c \in \Hom(\bQ_{l}, H^{2p}(X/k, \bQ_{l})(p)),
$$
and let
$ c' \in \Hom(\bQ_{l}, H^{2p}(\tX/k, \bQ_{l})(p)) $ be its
pull-back to
$ \tX $.
It is enough to show
$ c' $ algebraic.
Using the Leray spectral sequence and passing to the generic point of
$ S $,
$ c' $ induces
$$
c'' \in \Hom(\bQ_{l}, H^{2p}(Y/K, \bQ_{l})(p)).
$$
By
$ \TC(Y/K,p) $,
we get a cycle on
$ Y $ whose cycle class is
$ c'' $.
This induces a cycle on
$ \tX $ whose cycle class coincides with
$ c' $ mod
$ L^{0} $.
Here
$ c' $ mod
$ L^{0} $ belongs to
$$
\aligned
&\Hom(\bQ_{l}, H^{-1}(S/k, H^{2p-n+1}f_{*}(\bQ_{l,\tX}[n])(p)))
\\
&=
\Hom(\bQ_{l,S}[1], H^{2p-n+1}f_{*}(\bQ_{l,\tX}[n])(p)).
\endaligned
$$
So we may assume
$ c' \in L^{0} $ by modifying
$ c' $.
Then
$ c' \in L^{1} $ by (2.4.2), and the assertion follows from
$ \TC(X_{s}/k(s),p-1) $ using (2.4.4).
Here we may assume
$ k(s) = k $ (after reducing to the case
$ c' \in L^{1}) $ by the same argument as above.

\medskip\noindent
{\bf 2.6.~Proof of (0.7).}
Let
$ c \in \Hom(\bQ_{l}, H^{n}(X/k, \bQ_{l})(p)) $,
and
$ c' $ as above.
It is enough to show
$ c' $ algebraic.
We first reduce to the case
$ c' \in L^{0} $.
For this, we may assume
$ k(s) = k $ with
$ s $ as in condition~(i) (replacing
$ X $ with
$ X\otimes_{k}k(s) $,
and using an argument similar to Remark~(2.2)(ii)), because the
Leray spectral sequence is compatible with the base extension by
$ k \to k' $.
Let
$ i : X_{s} \to X $ denote the inclusion morphism.
Then
$ c' \in L^{0} $ is equivalent to
$$
c \in \Ker(i^{*} : H^{n}(X/k, \bQ_{l})(p) \to H^{n}(X_{s}/k,
\bQ_{l})(p)),
$$
by (2.4.3).
Using condition~(i) together with the hard Lefschetz theorem, we may
assume
$ c' \in L^{0} $ by modifying
$ c' $, because
$$
i^{*}i_{*} : H^{n-2}(X_{s}/k, \bQ_{l})\to H^{n}(X_{s}/k,
\bQ_{l})(1)
$$
coincides with the cup product with the hyperplane section class.
(The last fact can be verified by taking the composition with
$ i^{*} : H^{n-2}(X/k, \bQ_{l})\to H^{n-2}(X_{s}/k, \bQ_{l}) $
and
$ i_{*} : H^{n}(X_s/k, \bQ_{l})(1)\to H^{n+2}(X/k, \bQ_{l})(2) $
which are isomorphisms.)

We now reduce the assertion to the case
$ c' \in L^{1} $.
By the Poincar\'e duality, we have a canonical pairing on the
$ l $-adic sheaf
$ H^{0}f_{*}(\bQ_{l,\tX}[n])_{\ok}[-1] $ (see (2.4.1)),
which corresponds to a self duality
isomorphism of
$ H^{0}f_{*}(\bQ_{l,\tX}[n]) $.
Using this, we get a canonical decomposition in
$ \Perv(S/k,\bQ_{l}) $:
$$
H^{0}f_{*}(\bQ_{l,\tX}[n]) = \cF_{\inv} \oplus \cF_{\van}
$$
such that
$ (\cF_{\inv})_{\ok}[-1] $ is a constant sheaf and
$ \Gamma (S_{\ok}, (\cF_{\van})_{\ok}[-1]) = 0 $,
i.e.,
$$
H^{-1}(S/k,\cF_{\van}) = 0.
$$
Then
$ H^{1}(S/k,\cF_{\van}) = 0 $ by the duality.
Since
$ S = \bP^{1} $,
we have
$$
H^{0}(S/k,\cF_{\inv}) = 0,
$$
(i.e.,
$ H^{1}(S_{\ok}, (\cF_{\inv})_{\ok}[-1]) = 0) $,
and we get
$$
H^{0}(S/k, H^{0}f_{*}(\bQ_{l,\tX}[n])) = H^{0}(S/k,\cF_{\van}).
$$

So
$ c' $ mod
$ L^{1} $ is uniquely lifted to
$$
c'' \in \Hom(\bQ_{l}, (a_{S})_{*}\cF_{\van}),
$$
which corresponds to
$$
e \in \Ext^{1}(\bQ_{l,S}[1], \cF_{\van}),
$$
by the adjunction isomorphism for
$ a_{S} : S \to \Spec k $.
Using the restriction morphism by
$ Y \to X_{K} $,
we define
$ H^{j}(Y/K, \bQ_{l})_{X_{K}} $,
$ H^{j}(X_{K}/K, \bQ_{l})^{Y} $ as in Introduction so that we
have an exact sequence
$$
0 \to H^{n-1}(Y/K, \bQ_{l})_{X_{K}} \to
{H}_{c}^{n}((X_{K} \backslash Y)/K, \bQ_{l})\to
H^{n}(X_{K}/K, \bQ_{l})^{Y} \to 0.
$$
Then
$ H^{n-1}(Y/K, \bQ_{l})_{X_{K}} $ is isomorphic to the restriction of
$ \cF_{\van}[-1] $ to the generic point
$ \Spec K $ of
$ S $ (using (1.16)), and the restriction of
$ e $ to
$ \Spec K $ coincides with the pull-back of the above short exact sequence
by
$ c $.

Indeed, the adjunction isomorphism for
$ a_{S} : S \to \Spec k $ is induced by the morphism
$$
(a_{S})^{*}(a_{S})_{*} = (p_{2})_{*}(p_{1})^{*} \to
(p_{2})_{*}\delta_{*}\delta^{*}(p_{1})^{*} = id
$$
where
$ p_{a} : S\times_{k}S \to $
$ S $ are canonical projections, and
$ \delta : S \to S\times_{k}S $ is a diagonal embedding.
We get the coincidence using the exact sequence as above with
$ X $ replaced by
$ \tX $ together with a canonical morphism of the exact sequences
induced by the restriction morphism for
$ \pi : \tX \to X $.

So condition~(ii) (together with (1.16)) implies that there exist a
nonempty
open subvariety
$ S' $ of
$ S $ and a cycle
$ \zeta $ on
$ f^{-1}(S') $ such that the cycle class of the restriction of
$ \zeta $ to the generic fiber is zero and the restriction of
$ c' $ mod
$ L^{1} $ to
$ f^{-1}(S') $ is induced by
$ \zeta $.
Since
$ H^{0}f_{*}(\bQ_{l,\tX}[n]) $ is an intersection complex $\cF$ with
strict support
$ S $,
the natural morphism
$$
H^{0}(S/k, H^{0}f_{*}(\bQ_{l,\tX}[n])) \to H^{0}(S'/k,
H^{0}f_{*}(\bQ_{l,\tX}[n]))
$$
is injective.
(Indeed, $j_*j^*\cF/\cF$ is supported on $S\setminus S'$
where $j:S'\to S$ is the inclusion, and hence its cohomology
vanishes except for the degree 0).
So we may assume
$ S' = S $ by taking an extension of
$ \zeta $ to
$ \tX $,
and the assertion is reduced to the case
$ c' \in L^{1} $ by modifying
$ c' $.
Then we may assume again
$ k(s) = k $ by the same argument as in Remark~(2.2)(ii), and the
assertion follows from condition~(i) and (2.4.4).

\bigskip\bigskip
\centerline{{\bf 3. Divisor Case}}

\bigskip\noindent
In this section we treat the divisor case, and relate the
Tate conjecture to the de Rham conjecture for nonproper varieties
and finiteness of the Tate-Shafarevich groups.

\medskip\noindent
{\bf 3.1.}
Let
$ k $ be a finitely generated field over
$ \bQ $,
$ S $ an integral curve over
$ k $,
and
$ A $ an abelian scheme over a smooth dense open subvariety
$ U $ of
$ S $.
Let
$ \cL = T_{l}A $ which is the projective system of \'etale sheaves defined
by
$ \Ker(l^{m} : A \to A) $.
We assume
$$
H^{0}(U_{\ok}, \cL) = 0.
\leqno(3.1.1)
$$

Let
$ K = k(S) $,
and
$ \tK $ be the maximal subfield of
$ \oK $ that is unramified over
$ U $.
Let
$ G = \Gal(\tK/K) $ which is identified with
$ \pi_{1}(U, \Spec \oK) $.
Then
$ \cL $ corresponds to the
$ \bZ_{l} $-module
$$
E := \Gamma (\Spec \tK, \cL) = T_{l}A(\tK)
$$
with a continuous action of
$ G $,
and we have natural isomorphisms
$$
H^{1}(G, E/l^{m}) = H^{1}(U,\cL/l^{m}),\quad H^{1}(G, E) = H^{1}(U,\cL).
\leqno(3.1.2)
$$
Here
$ \cL/l^{m} $ denotes
$ \cL/l^{m}\cL $ for simplicity (same for
$ E) $,
and
$ H^{1}(G, E) $ is the projective limit of
$ H^{1}(G, E/l^{m}) $ (same for
$ H^{1}(U,\cL)) $.
Using the Kummer sequence, we get an exact sequence
$$
0 \to A(K)\otimes \bZ_{l}\to H^{1}(G, E) \to
T_{l}H^{1}(G, A(\tK)) \to 0.
\leqno(3.1.3)
$$
Note that
$ A(K) $ is a finitely generated abelian group by (a generalization of)
the Mordell-Weil theorem, and
$ H^{1}(G, E) $ is a finite
$ \bZ_{l} $-module by (3.1.5) below.

Let
$ m $ be a positive integer.
By the Hochschild-Serre spectral sequence, we have a long exact sequence
$$
\aligned
0
&\to H^{1}(G_{k}, H^{0}(U_{\ok}, \cL/l^{m})) \to
H^{1}(U, \cL/l^{m}) \to H^{0}(G_{k}, H^{1}(U_{\ok}, \cL/l^{m}))
\\
&\to H^{2}(G_{k}, H^{0}(U_{\ok}, \cL/l^{m})).
\endaligned
\leqno(3.1.4)
$$
Passing to the limit, we get an exact sequence by the Mittag-Leffler
condition, and
$$
H^{1}(U,\cL) = H^{1}(U_{\ok},\cL)^{G_{k}},
\leqno(3.1.5)
$$
because (3.1.1) implies that
$ H^{0}(U_{\ok}, \cL/l^{m+k}) \to H^{0}(U_{\ok}, \cL/l^{m}) $
is zero for
$ k \gg 0 $.

\medskip\noindent
{\bf 3.2.~Remarks.} (i) There exists a positive integer
$ a $ independent of
$ m $ such that
$ H^{0}(U_{\ok}, \cL/l^{m}) $ is annihilated by
$ l^{a} $,
because
$ H^{0}(U_{\ok}, \cL\otimes_{\bZ_{l}}(\bQ_{l}/\bZ_{l}))
$ is finite by (3.1.1) (using also the finiteness of
$ H^{1}(U_{\ok}, \cL) $ over
$ \bZ_{l}) $.

\medskip
(ii) Assume
$ S = \bP^{1} $.
Let
$ b $ be a positive integer such that there is no elements of
$ A(K) $ with order
$ l^{b} $.
Then there exists a thin subset
$ \Sigma $ of
$ U(k) $ in the sense of [33] such that there is no element of
$ A_{s}(k) $ with order
$ l^{b} $ for
$ s \in U(k) \backslash \Sigma $,
where
$ A_{s} $ is the fiber of
$ A $ over
$ s $.
(It is enough to apply Hilbert's irreducibility theorem ([21], [33])
to the irreducible components of the kernel of
$ l^{b} : A \to A $.)
In particular, the
$ l $-primary torsion part of
$ A_{s}(k) $ is annihilated by
$ l^{b-1} $.
We have the same for the torsion part of
$ H^{1}(G_{s}, E) $ using the exact sequence
$$
0 \to A_{s}(k)\otimes \bZ_{l}\to H^{1}(G_{s}, E)
\to T_{l}H^{1}(G_{s}, A(\ok)) \to 0,
\leqno(3.2.1)
$$
because the last term is torsion-free.

\medskip
(iii) The injectivity of (0.10) in Introduction follows immediately
from
Hilbert's irreducibility theorem (see [21], [33]), if
$ V $ is replaced with
$ E/l^{m} $, see also (3.3.2) below.
Taking the projective limit for $m$, we have the assertion for
$ E $.
(For this we do not have to assume that
$ E $ is associated with an abelian variety.)
Then the injectivity of (0.10) for
$ V $ (with
$ |U| $ replaced by
$ U(k)\setminus\Sigma$) follows from Remark~(ii) above since we may assume that the
thin subsets used in Hilbert's irreducibility theory for the above
argument contain the thin subset $\Sigma$ in Remark~(ii) above.
Note that the assertion is true without assuming
$ S = \bP^{1} $ (using a morphism to
$ \bP^{1}) $.
A.~Tamagawa has informed us of the injectivity for any smooth sheaf
$ \cL $ (not necessarily associated with an abelian scheme) using the
theory of Frattini subgroups, see also [30].
We are also informed that a similar argument was used in [38].

\medskip\noindent
{\bf 3.3.~Proposition.} {\it
With the notation and the assumptions of {\rm (3.1)}, assume further
$ S = \bP^{1} $.
Then there exists a thin subset
$ \Sigma $ of
$ U(k) $ in the sense of {\rm [33]} such that
$ H^{1}(G,E) \to H^{1}(G_{s},E) $ is injective for
$ s \in U(k) \backslash \Sigma  $, where
$ G_{s} = \Gal(\overline{k(s)}/k(s)) $.
}

\medskip\noindent
{\it Proof.}
Let
$ a, b $ be positive integers as in Remarks~(i) and (ii) above,
and take an integer
$ m > a + b $.
For a subfield
$ K' $ of
$ \tK $ which is Galois over
$ K $,
let
$ U' \to U $ be the corresponding finite \'etale morphism.
For each
$ s \in |U| $,
we choose a lift
$ s' \in |U'| $ in a compatible way with the change of
$ K' $ so that we get an inclusion
$ G_{s} \to G $ and a morphism of (3.1.3) to (3.2.1).
We assume
$ K' $ sufficiently large so that the pull-back of
$ \cL/l^{m} $ to
$ U' $ is trivial, and the image of the composition of
$$
\aligned
&H^{1}(\Gal(K'/K), E/l^{m})
\to H^{1}(G, E/l^{m})
\\
&\quad = H^{1}(U,\cL/l^{m})\to H^{1}(U_{\ok}, \cL/l^{m})^{G_{k}}
\endaligned
\leqno(3.3.1)
$$
coincides with the image of
$ H^{1}(G, E/l^{m}) \to H^{1}(U_{\ok}, \cL/l^{m})^{G_{k}} $.
Combined with (3.1.4) and Remark~(i) above, this implies that for any
$ e \in H^{1}(G, E/l^{m}) $, there exists
$ e' \in H^{1}(\Gal(K'/K), E/l^{m}) $ such that
$ e - e' $ is annihilated by
$ l^{a} $.
(Since the first morphism of (3.3.1) is injective, the image of
$ e' $ is also denoted by
$ e' $.)

By Hilbert's irreducibility theorem, there exists a thin subset
$ \Sigma $ of
$ U(k) $ such that for
$ s \in U(k) \backslash \Sigma $,
a point
$ s' $ of
$ U' $ over
$ s $ is unique and hence
$ \Gal(K'/K) = \Gal(k(s')/k(s)) $, see [21], [33].
So we get the isomorphism
$$
H^{1}(\Gal(K'/K), E/l^{m}) \simto H^{1}(\Gal(k(s')/k(s)), E/l^{m}),
\leqno(3.3.2)
$$
and hence any
$ e\in\Ker(H^{1}(G, E/l^{m}) \to H^{1}(G_{s}, E/l^{m})) $ is
annihilated by
$ l^{a} $ (taking $e'$ as above).
We may assume moreover the torsion part of
$ H^{1}(G_{s}, E) $ is annihilated by
$ l^{b-1} $ for
$ s \in U(k) \backslash \Sigma $ by Remark~(ii) above, replacing
$ \Sigma $ if necessary.
Set
$$
H = H^{1}(G, E),\quad H' = H^{1}(G_{s}, E).
$$
Then
$ \Ker(H/l^{m} \to H'/l^{m}) $ is also annihilated by
$ l^{a} $, because the natural morphism
$$
H^{1}(G, E)/l^{m} \to H^{1}(G, E/l^{m})
$$
is injective.
(This can be verified easily by taking the limit of the exact sequence
associated
with
$ 0 \to E/l^{i} \to E/l^{i+m} \to E/l^{m}
\to 0 $.)
Since
$ \Ker(A(K) \to A_{s}(k)) $ and
$ T_{l}H^{1}(G, A(\tK)) $ are torsion-free,
$ \Ker(H \to H') $ is also torsion-free by (3.1.3) and (3.2.1).
Finally,
$ H $ is finite over
$ \bZ_{l} $ by (3.1.2) and (3.1.5).
So the assertion follows from

\medskip\noindent
{\bf Sublemma.}
A morphism of
$ \bZ_{l} $-modules
$ H \to H' $ is injective, if
$ \Ker(H \to H') $ is torsion-free,
$ H $ is finite over
$ \bZ_{l} $,
and
$ \Ker(H/l^{m} \to H'/l^{m}) $ and the torsion part of
$ H' $ are annihilated by
$ l^{a} $ and
$ l^{b} $ respectively with
$ a + b < m $.

\medskip\noindent
{\it Proof.}
We may assume the surjectivity of the morphism by replacing
$ H' $ with the image so that we get an exact sequence. (Here
$ \Ker(H/l^{m} \to H'/l^{m}) $ becomes smaller by replacing $H'$.)
Then the assertion follows from the snake lemma applied to the
endomorphism of the short exact sequence defined by the
multiplication by
$ l^{m} $.

\medskip\noindent
{\bf 3.4.~Remarks.} (i)
Let
$ f : X \to S $ be a surjective morphism of smooth irreducible
projective
varieties over
$ k $.
(We may assume
$ \dim X = 2 $ since the Tate conjecture for divisors is
reduced to the surface case by (0.5) replacing
$ k $ if necessary.)
Let
$ U $ be a dense open subvariety of
$ S $ such that the restriction
$ f' : X' \to U $ of
$ f $ over
$ U $ is smooth.
We have a closed embedding
$ X' \to X\times U $ by the graph of
$ f $.
Let
$$
\aligned
A &= \Coker(J_{X\times U/U} \to J_{X'/U}),
\\
\cL &= \Coker(R^{1}(pr_{2})_{*}\bZ_{l,X\times U}(1) \to
R^{1}f'_{*}\bZ_{l,X'}(1)),
\endaligned
$$
where
$ J_{X\times U/U}, J_{X'/U} $ denote the identify component of
the Picard scheme, and the cokernel
$ A $ can be defined by using the N\'eron model of the cokernel over the
generic point of
$ U $.
Then we have a canonical isomorphism
$ \cL = T_{l}A $,
and (3.1.1) is satisfied by the global invariant cycle theorem.

\medskip
(ii) With the notation of Remark~(i) above, the Tate conjecture for
divisors on
$ X $ is equivalent to
$$
T_{l}H^{1}(G, A(\tK)) = 0.
\leqno(3.4.1)
$$
Indeed, let
$ Y $ be the generic fiber of
$ X \to S $,
and
$ \CH_{\hom}^{1}(Y) $ as in Introduction.
Since the Leray spectral sequence for
$ X'_{\ok} \to U_{\ok} \to \Spec
\ok $ degenerates at
$ E_{2} $ (see [6]), the cycle map induces the surjective morphism
$$
\CH_{\hom}^{1}(Y)\otimes \bQ_{l}\to
H^{1}(U_{\ok},\cL)^{G_{k}}\otimes_{\bZ_{l}}\bQ_{l},
$$
if the Tate conjecture is true.
It is known that this morphism coincides with the composition of the
canonical morphism
$ \CH_{\hom}^{1}(Y) \to A(K) $ with the first morphism of (3.1.3)
with
$ \bQ_{l} $-coefficients.
So we get the assertion by (3.1.3), because
$ T_{l}H^{1}(G, A(\tK)) $ is torsion-free.

\medskip\noindent
{\bf 3.5.~Remark} (S.~Bloch, K.~Kato, J.~Nekovar).
Let
$ k $ be a number field, and
$ A $ an abelian variety over
$ k $.
Let
$ E = T_{l}A(\ok) $.
For a place
$ v $ of
$ k $,
let
$ k_{v} $ be the completion of
$ k $ at
$ v $,
$ \ok_{v} $ its algebraic closure, and
$ E_{v} = T_{l}A(\ok_{v}) $,
where
$ E_{v} $ is isomorphic to
$ E $ by choosing
$ \ok \to \ok_{v} $.
Then we have a canonical morphism of exact sequences
$$
\CD
0 @>>> A(k)\otimes \bZ_{l} @>>> H^{1}(G_{k}, E) @>>>
T_{l}H^{1}(G_{k}, A(\ok)) @>>> 0
\\
@. @VVV @VVV @VVV
\\
0 @>>>A(k_{v})\hat{\otimes }\bZ_{l} @>>> H^{1}(G_{k_{v}}, E_{v})
@>>> T_{l}H^{1}(G_{k_{v}}, A(\ok_{v})) @>>> 0.
\endCD
$$
If
$ e \in H^{1}(G_{k}, E) $ is defined geometrically by taking a pull-back
of
an exact sequence like (0.2) in Introduction (where
$ A = J_{D,X}) $,
then we have the same for the restriction
$ e_{v} \in H^{1}(G_{k_{v}}, E_{v}) $ of
$ e $,
and
$ e_{v} $ belongs to the image of
$ A(k_{v})\hat{\otimes }\bZ_{l} $ for any
$ v $,
using the theory of Bloch and Kato  on
$ {H}_{g}^{1} $ (see [3], Remark before 3.8),
provided that the de Rham conjecture [15] holds for the nonproper
variety
$ (X \backslash D)_{v} $, see also [24].
If this is true, then
$ e $ belongs to the image of
$ A(k)\otimes \bZ_{l} $ if we have furthermore the injectivity of
$$
T_{l}H^{1}(G_{k}, A(\ok)) \to \prod_{v} T_{l}H^{1}(G_{k_{v}},
A(\ok_{v})),
$$
(i.e., the
$ l $-primary torsion part of the Tate-Shafarevich group of
$ A $ is finite).


\begin{thebibliography}{99}

\bibitem{[1]}
Beilinson, A., On the derived category of perverse sheaves, Lect. Notes
in Math., vol. 1289, Springer, Berlin, 1987,. pp. 27--41

\bibitem{[2]}
Beilinson, A., Bernstein, J. and Deligne, P., Faisceaux pervers,
Ast\'erisque, vol. 100, Soc. Math. France, Paris, 1982.

\bibitem{[3]}
Bloch, S. and Kato, K., L-functions and Tamagawa numbers of motives,
in Grothendieck Festschrift I, Birkh\"auser, Boston, 1990,
pp. 333--400.

\bibitem{[4]}
Cassels, J. and Schinzel, A., Selmer's conjecture and families of
elliptic curves, Bull. London Math. Soc. 14 (1982), 345--348.

\bibitem{[5]}
Conrad, B., Conrad, K. and Helfgott, H.,
Root numbers and ranks in positive characteristic,
Adv. Math. 198 (2005), 684--731.

\bibitem{[6]}
Deligne, P., Th\'eor\`eme de Lefschetz et crit\`eres
de d\'eg\'en\'erescence de suites spectrales, Publ. Math. IHES, 35
(1968), 107--126.

\bibitem{[7]}
Deligne, P., Conjecture de Weil II, Publ. Math. IHES, 52 (1980),
137--252.

\bibitem{[8]}
Deligne, P., Le groupe fondamental de la droite projective moins
trois points, in Galois groups over Q, Springer, New York, 1989, pp.
79--297.

\bibitem{[9]}
Deligne, P., La formule de Picard-Lefschetz, in Lect. Notes in Math.,
vol. 340, Springer Berlin, 1973, pp. 165--196.

\bibitem{[10]}
Deligne, P., Le formalisme des cycles \'evanescents,
in SGA7 XIII and XIV, Lect. Notes in Math., vol. 340, Springer, Berlin,
1973, pp. 82--115 and 116--164.

\bibitem{[11]}
Deligne, P., Th\'eor\`eme de finitude en cohomologie {\it l}-adique,
in SGA 4 1/2, Lect. Notes in Math., vol. 569, Springer, Berlin, 1977,
233--261.

\bibitem{[12]}
Deligne, P., Milne, J., Ogus, A. and Shih, K., Hodge cycles, motives, and
Shimura varieties, Lect. Notes in Math., vol. 900, Springer, Berlin, 1982.

\bibitem{[13]}
Ekedahl, T., On the adic formalism, in Grothendieck Festschrift II,
Birkh\"auser, Boston, 1990, pp. 197--218.

\bibitem{[14]}
El Zein, F and Zucker, S.,
Extendability of normal functions associated to algebraic
cycles, in Topics in transcendental algebraic geometry, Ann. Math.
Stud., 106, Princeton Univ. Press, Princeton, N.J., 1984, pp. 269--288.

\bibitem{[15]}
Faltings, G., Crystalline cohomology and
$ p $-adic Galois representation, in Algebraic Analysis, Geometry and
Number Theory, edited by J. Igusa, Johns Hopkins University Press,
Baltimore, 1989, pp. 25--80.

\bibitem{[16]}
Griffiths, P., On the period of certain rational integrals I,
II, Ann. Math. 90 (1969), 460--541.

\bibitem{[17]}
Honda, T., On the theory of formal groups, J. Math. Soc. Japan 22 (1970),
213--246.

\bibitem{[18]}
Houzel, C., Morphisme de Frobenius et rationalit\'e des fonctions
$ L $, in SGA 5, Lect. Notes in Math., vol. 589, Springer, Berlin, 1977,
pp. 442--480.

\bibitem{[19]}
Jannsen, U., Mixed motives and algebraic K-theory, Lect. Notes in Math.,
vol. 1400, Springer, Berlin, 1990.

\bibitem{[20]}
Katz, N., Etude cohomologique des pinceaux de Lefschetz, in Lect. Notes
in Math., vol. 340, Springer Berlin, 1973, pp. 254--327.

\bibitem{[21]}
Lang, S., Fundamentals of Diophantine geometry, Springer, Berlin,
1983.

\bibitem{[22]}
Lefschetz, S., L'analysis situs et la g\'eom\'etrie
alg\'ebrique, Paris Gauthier-Villars, 1924.

\bibitem{[23]}
Murre, J.P., Lectures on an Introduction to Grothendieck's Theory of
the Fundamental Group, Lecture Notes, Tata Institute of Fundamental
Research, Bombay, 1967.

\bibitem{[24]}
Nekov\'a\v r, J.,
$p$-adic Abel-Jacobi maps and $p$-adic heights,
in The arithmetic and geometry of algebraic cycles (Banff, AB, 1998),
CRM Proc. Lecture Notes, 24, Amer. Math. Soc., Providence, RI, 2000,
pp. 367--379.

\bibitem{[25]}
Saito, M., Mixed Hodge Modules, Publ. RIMS Kyoto Univ. 26 (1990),
221--333.

\bibitem{[26]}
Saito, M., Hodge conjecture and mixed motives, I, Proc. Symp. Pure
Math. 53 (1991), 283--303.

\bibitem{[27]}
Saito, M., On the formalism of mixed sheaves, preprint RIMS-784,
Aug. 1991 (also available at arXiv:math/0611597).

\bibitem{[28]}
Saito, M., Admissible normal functions, J. Alg. Geom., 5 (1996),
235--276.

\bibitem{[29]}
Selmer, E., A conjecture concerning rational points on cubic curves,
Math. Scand., 2 (1954), 49--54.

\bibitem{[30]}
Serre, J.-P., Groupes analytiques
$ p $-adiques, S\'eminaire Bourbaki, 270 (1963/64).

\bibitem{[31]}
Serre, J.-P., Abelian {\it l}-adic representations and elliptic curves,
Benjamin, New York, 1968.

\bibitem{[32]}
Serre, J.-P., Groupes
$ p $-divisibles (d'apr\`es J. Tate), S\'eminaire Bourbaki 1966/67 no 318.

\bibitem{[33]}
Serre, J.-P., Lectures on the Mordell-Weil theorem, Vieweg,
Braunschweig, 1989.

\bibitem{[34]}
Steenbrink, J. and Zucker, S., Variation of mixed Hodge structure I, Inv.
Math., 80 (1985), 489---542.

\bibitem{[35]}
Tate, J., Algebraic cycles and poles of zeta functions, in Arithmetic
algebraic geometry, Harper \& Row, New York, 1965, pp. 93--110.

\bibitem{[36]}
Tate, J., On the conjecture of Birch and Swinnerton-Dyner and a
geometric analog, S\'eminaire Bourbaki 1965/66 no 306 (reproduced
in Dix Expos\'es sur la cohomologie des schemas, North-Holland,
Amsterdam, 1968, pp. 289--214).

\bibitem{[37]}
Tate, J., Conjectures on algebraic cycles in $l$-adic cohomology,
Motives (Seattle, WA, 1991), Proc. Sympos. Pure Math., 55,
Part 1, Amer. Math. Soc., Providence, RI, 1994, pp. 71--83.

\bibitem{[38]}
Terasoma, T., Complete intersections with middle Picard number
$ 1 $ defined over {\bf Q}, Math. Z. 189 (1985), 289--296.

\bibitem{[39]}
Verdier, J.-L., Extension of a perverse sheaf over a closed subspace,
Ast\'erisque, 130 (1985), 210--217.

\bibitem{[40]}
Verdier, J.-L., Cat\'egories d\'eriv\'ees, in SGA 4 1/2, Lect. Notes
in Math., vol. 569, Springer, Berlin, 1977, 262--311.

\bibitem{[41]}
Zucker, S., The Hodge conjecture for cubic fourfolds, Compos. Math., 34
(1977), 199--209.

\end{thebibliography}
\end{document}